\documentclass[12pt]{amsart}
\usepackage{fullpage} 
\usepackage{amssymb}
\usepackage{amsmath}
\usepackage[dvipsnames]{xcolor}
\usepackage{graphicx,url}        % standard LaTeX graphics tool
\newcommand{\tk}{\tilde{k}}
\newcommand{\tc}{\tilde{c}}
\newcommand{\ttk}{\tilde{\tilde{k}}}
\newcommand{\ttc}{\tilde{\tilde{c}}}
\newcommand{\Um}{U^{(m)}}
\newcommand{\Umm}{U^{(m-1)}}
\newcommand{\Ummm}{U^{(m-2)}}
\newcommand{\bbar}{\bar{b}}

\newcommand{\vz}{v_{init}}
\newcommand{\uz}{u_{init}}
\newcommand{\ddt}{\frac{d}{dt}}

\newcommand{\Chi}{\Upsilon}

\newcommand{\mpcitee}[2]{{\cite{#1}(#2)}}
\newcommand{\fF}{\mathcal{F}}

\newcommand{\cC}{\mathcal{C}}

\newcommand{\inner}[1]{\langle#1\rangle}

\newcommand{\bsmult}{\begin{multline}}
\newcommand{\esmult}{\end{multline}}
\newcommand{\QT}{Q^{\mathcal{T}}}

\newcommand{\LambdaT}{\Lambda^{\mathcal{T}}}
\newcommand{\LfF}{L_{\fF}}
\newcommand{\lambdaT}{\lambda^{\mathcal{T}}}
\newcommand{\hup}{\gamma_{r}}
\newcommand{\hdown}{\gamma_{l}}
\newcommand{\fl}{\hdown}
\newcommand{\fr}{\hup}

\newcommand{\mpc}[1]{{#1}}
\newcommand{\mpa}[1]{{#1}}

\newcommand{\bsa}{\begin{subequations}}
\newcommand{\esa}{\end{subequations}}
\newtheorem{example}{Example}%[section]
\newtheorem{theorem}{Theorem}%[section]
\newtheorem{prop}{Proposition}%[section]
\newtheorem{lemma}{Lemma}%[section]
\newtheorem{remark}{Remark}%[section]
\newtheorem{assum}{Assumption}%[section]
\newtheorem{corollary}{Corollary}%[section]
\newcommand{\rthe}[1]{Theorem~\ref{#1}}

\newcommand{\rlem}[1]{Lemma~\ref{#1}}

\newcommand{\ba}{\begin{eqnarray}}
\newcommand{\ea}{\end{eqnarray}}
\newcommand{\bsub}{\begin{subequations}}
\newcommand{\esub}{\end{subequations}}
\newcommand{\bas}{\begin{eqnarray*}}
\newcommand{\eas}{\end{eqnarray*}}
\newcommand{\mpunit}[1]{\mathrm{[#1]}}
\newcommand{\ourcite}[2]{\cite{#1}~(#2)}
\newcommand{\tT}{\mathcal{T}}

\newcommand{\cg}{\mathcal{C}}
\newcommand{\A}{\mathbb{A}}

\newcommand{\vbar}{\overline{V}}

\newcommand{\Hi}{V}

\newcommand{\cab}[1]{\cg(a,b;#1)}
\newcommand{\ccab}{\cg(a,b)}
\newcommand{\rab}[1]{R(a,b;#1)}
\newcommand{\resab}[3]{R(#1,#2;\;#3)}

\newcommand{\Dom}{\operatorname{Dom}}
\newcommand{\Rg}{\operatorname{Rg}}

\newcommand\myskip[1]{}

\newcommand{\R}{\mathbb{R}}
\newcommand{\abs}[1]{\mid\!#1\!\mid}
\newcommand{\norm}[2]{||#1||_{#2}}
\newcommand{\jJ}{\mathcal{J}}

\title[Algorithms for non-equilibrium and hysteresis in permafrost models]
{
Unified analysis of algorithms\\
for equilibrium, non-equilibrium, and hysteresis models of phase transition in permafrost}
\author{Malgorzata Peszynska}
\author{Nicholas Slugg}
\address{Oregon State University, Department of Mathematics, Corvallis, OR 97331, USA} \email{mpesz@math.oregonstate.edu, sluggn@oregonstate.edu}
\begin{document}
\begin{abstract}
In this paper we consider a nonlinear partial differential equation describing heat flow with ice-water phase transition in permafrost soils. Such models and their numerical approximations have been well explored in the applications literature. In this paper we describe  a new direction in which the allow relaxation and hysteresis of the phase transition which introduce additional nonlinear terms and complications for the analysis. We present numerical algorithms  as well as analysis of the well-posedness and convergence of the fully implicit iterative schemes. The analysis we propose handles the equilibrium, non-equilibrium, and hysteresis cases in a unified way.
 We also illustrate with numerical examples for a model ODE and PDE. 
 \end{abstract}
 \maketitle

%%%%%%%%%%%%%%%
\section{Introduction}

In this paper we study algorithms for approximation of solutions to a nonlinear heat equation modeling the thermal processes of thawing and freezing in a soil. Such processes are important for modeling in permafrost regions (e.g., in the Arctic) \mpc{which} respond to the daily and seasonal variation of temperature at the ground surface.  
The specific challenge we consider in this paper is that the phase change \mpc{processes} (from liquid to ice and vice-versa) need not be instantaneous. Rather, they are assumed to follow non-equilibrium (also known as kinetic or relaxation laws) or hysteretic relationships. 

The models we consider are nonlinear, and the appropriate numerical schemes for the  approximation of solutions already for the equilibrium case present interesting challenges well documented in our previous work \cite{BPV,VP-Tp}.    
In the paper we propose numerical discretization algorithms extending those for the equilibrium relationships to the cases of non-equilibrium and hysteretic relationships. We also analyze them within a common framework. In particular, we provide results on the well-posedness of the discrete schemes as well as conditions on the convergence of iterative algorithms used to solve the nonlinear problem. 

Specifically, we consider a family of nonlinear evolution problems in \mpa{an open bounded domain} $\Omega \subset \R^d$
\ba
\label{eq:pdeu}
\partial_t (c(u) + \chi) - \nabla \cdot (k(u) \nabla  u) = f(x,t), \; x \in \Omega, t>0.
\ea
The unknown $u$ in \eqref{eq:uchi} represents the temperature of thawing or freezing soil, i.e, undergoing phase transition. The term $c(u)$ is associated with variable heat capacity, and the coefficient $k(u)$ represents heat conductivity, and is a symmetric uniformly positive definite diffusivity coefficient. One also needs  initial conditions \mpa{at $t=0$ imposed on $c(u)+\chi$} and boundary conditions \mpa{imposed on the boundary $\partial \Omega$ of $\Omega$ which is assumed to be Lipschitz~\cite{Visbook}}. 

To complete \eqref{eq:pdeu} we need to define the relationship between $\chi$ and $u$.  In the model \eqref{eq:pdeu} $\chi$ is associated with the phase transition and describes the heat amount necessary for thawing, so that $c(u)+\chi(u)$ represent  the total energy density. In this paper we consider exactly one of the relationships coupling $\chi$ to $u$ discussed below. 
\bsa
\label{eq:uchi}
\ba
\label{eq:EQ}
\mathrm{(EQ)}\;\;\; &\chi&= F(u),
\\
\label{eq:NEQ}
\mathrm{(KIN)}\;\;\; 
\partial_t &\chi& + B(\chi - F(u))=0,
\\
\label{eq:HYST}
\mathrm{(HYST)}\;\;\; 
\partial_t &\chi& + \cg(\chi -F(u)) \ni 0.
\ea
\esa
In these models $F$ is a bounded nonnegative function  which is monotone nondecreasing and smooth except at $u=0$ where it quickly grows to $1$. 

The choices in \eqref{eq:uchi} describe a phase transition that can occur instantaneously (in equilibrium) as in \eqref{eq:EQ}, or slowly. The non-equilibrium (kinetic) evolution model (KIN) \eqref{eq:NEQ} allows $\chi$ to evolve towards \eqref{eq:EQ} with some rate $B>0$. Such non-equilibrium models are common in the applications involving phase transition and chemical reactions; see our work in \cite{MedinaP18} and \cite{PShin21}. 
With the hysteretic model \eqref{eq:HYST} when $u$ increases, $\chi$ follows a different path than at freezing when $u$ decreases; this well-known phenomenon is due to the physics of nucleation \cite{Visbook,LittShow95}. The model represents evolution under constraints given by $\cg$, a multivalued monotone constraint graph.  

\begin{figure}[ht]
\centering
\includegraphics[width=.45\linewidth]{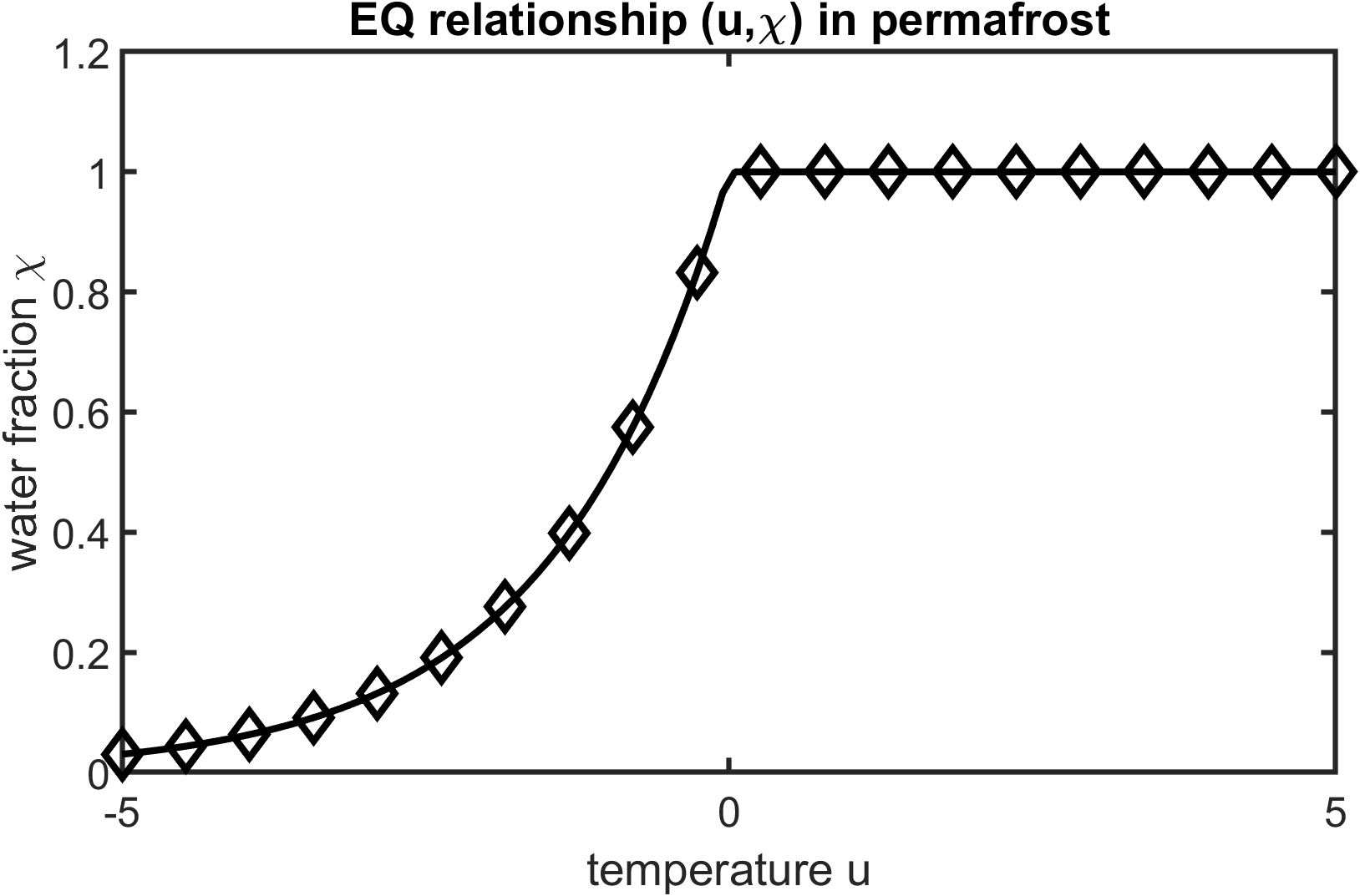}
\includegraphics[width=.45\linewidth]{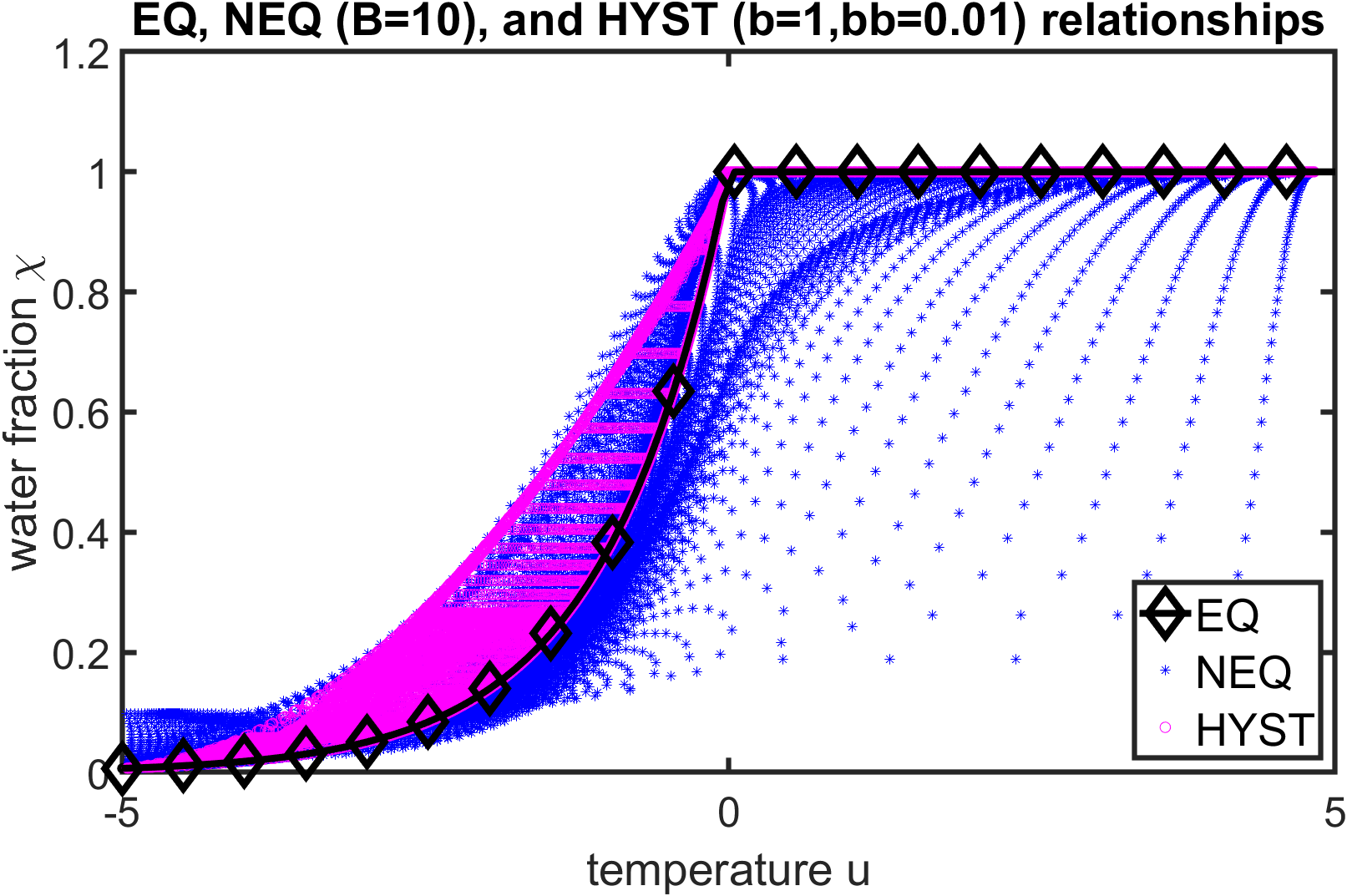}
%permafrost_noneq(0.1,30,0.1,20)
    \caption{Illustration of of relationships \eqref{eq:uchi}. On left, we plot of $\chi=F(u)$ in \eqref{eq:EQ}.  On right, we plot the same relationship $\chi=F(u)$ again as well as a scatter plot of the results $(u,\chi)$ of simulations+ for \eqref{eq:NEQ} and \eqref{eq:HYST}. Details of this example are given in Section~\ref{sec:results-pde}.}
    \label{fig:uchi-ALL}
\end{figure}

We aim to approximate the solutions $(u,\chi)$ satisfying  \eqref{eq:pdeu} and one of \eqref{eq:EQ}--\eqref{eq:HYST},  supplemented with appropriate boundary and initial conditions. For simplicity, we use finite differences (FD) in space, and in our examples we use specifically cell-centered FD \mpc{called} CCFD. In time, we use fully implicit first order scheme. The fully discrete model reads: at every time step $t_n$, solve
\ba
\label{eq:uchin}
\tfrac{1}{\tau}(C(U^n) -C(U^{n-1}))+ \tfrac{1}{\tau}(\Chi^n-\Chi^{n-1}) +A(U^n)U^n = f^n, 
\ea
where $U^n=(U_j^n)_j,\Chi^n=(\Chi^n_j)_j$ contain the approximations to $u(x_j,t_n),\chi(x_j,t_n)$ at spatial points $x_j \in \Omega$ of a uniform rectangular grid over $\Omega$. Here $C(U) =(c(U_j^n))_j$ is a vector, and \mpa{$A(U)$} is a \mpa{matrix which is symmetric positive definite for any $U$, and $A(U)U$ is an approximation to $-\nabla \cdot k(u) \nabla u$} under suitable boundary conditions, with $A(U)$ deriving its properties from $k=k(u)$. The system \eqref{eq:uchin} is complemented with appropriately discretized equations from \eqref{eq:uchi}. We refer to  \cite{BPV,VP-Tp} for details on $C,A$.

The resulting system is nonlinear and must be solved by some iteration, while its well-posedness is only guaranteed under some assumptions. The presence of multiple nonlinearities including $F(\cdot),c(\cdot),k(\cdot)$ as well as the complications due to non-equilibrium models makes this task challenging to implement and to analyze. The analysis in this paper addresses these difficulties in a unified way. We also present simulations for a few selected examples of \eqref{eq:uchi} and discuss robustness of the solver. 

We believe the results to be presented are new, even if they are related to some work in the literature on other applications.  In particular,  we are not aware of any computational models or rigorous work for permafrost applications with hysteresis. However, there is ample work on analysis and abstract framework on numerical schemes for the Stefan problem in equilibrium for Galerkin FEM and with relaxation, e.g., see \cite{Visbook,JiangNochetto98}. In our own prior work \cite{BPV,VP-Tp} we used CCFD but worked only with the equilibrium problem and considered only the time-lagging of \mpa{$A(U)$}.  

\medskip
{\bf Outline.} The outline of the paper is as follows. 
In Section~\ref{sec:prelim} we develop some notation and recall technical preliminaries.  In Section~\ref{sec:model} we provide details of the model \eqref{eq:uchi} and its numerical discretization. In Section~\ref{sec:analysis} we provide analysis of the well-posedness of the discrete schemes and of the iterative schemes.  In Section~\ref{sec:examples} we illustrate the results with examples for the simplified ODE counterpart of the model and for the PDE. In Section~\ref{sec:summary} we summarize and outline future and current work.

%%%%%%%%%%%%%%%%%
\section{Preliminaries and technical results} 
\label{sec:prelim}
In this Section we recall some material from the literature needed for the model developments and for the proofs to follow. In particular, we provide  information on the evolution problems under constraints.

%%%%
\subsection{Operator equation in a Hilbert space} 
 
In the paper we will consider a Hilbert space $\Hi$, with inner product $\inner{\cdot,\cdot}$, and norm $\norm{\cdot}{}$. When $\Hi=\R$, we will use $\inner{x,y}=xy,\norm{x}{}=\abs{x}$. On $\Hi=\R^n$ we will use the Euclidean inner product and norm.

Let $I$ be an identity operator on $\Hi$, $I(v)=v$.

For vector spaces $V,W$, and any function $T:V\to W$ we will denote by $L_T$ its Lipschitz constant:   
\ba
\norm{T(v_1)-T(v_2)}{W} \leq L_T \norm{v_1-v_2}{\mpc{V}}.
\ea
We say that $T$ is a contraction if $L_T<1$. 

Let $v_1,v_2 \in \Hi$ be arbitrary. We say that \mpc{$T:V\to V$} is monotone if
%%%
\ba
\inner{T(v_1)-T(v_2),v_1-v_2}\geq 0.
\ea
We also say that $T$ is strongly monotone if there is some $c_0>0$:
\ba
\inner{T(v_1)-T(v_2),v_1-v_2}\geq c_0\norm{v_1-v_2}{}^2.
\ea

The following equation is of interest in this paper
\ba
\label{eq:T}
T(u)=b.
\ea
%%%
Of interest is existence and uniqueness of the solutions to \eqref{eq:T}. Let $T:\Hi \to \Hi$. 

%%%%%%%%%
\begin{theorem}[\cite{AHan2}, Thm.~5.1.4, p211]
\label{th:AHAN}
Let $T:\Hi \to \Hi$ be (i) strongly monotone and (ii) Lipschitz continuous. Then for any $b\in \Hi$ there is a unique solution $u \in V$ to \eqref{eq:T}. 
Moreover, if $T(u_1)=b_1,T(u_2)=b_2$ then 
$
\norm{u_1-u_2}{} \leq \frac{1}{c_0}\norm{b_1-b_2}{}.
$.
\end{theorem}

The following result is also well known, and we will apply it in $K=V$. 
\begin{theorem}[\cite{AHan2}, Thm.~.5.1.3, p209]
\label{th:fixed}
Let $T:K\to K$ be a contraction where $K\subset V$ is a nonempty closed set. Then for any $v^0$ the sequence
\ba
v^m=T(v^{m-1}), m=1,2, \ldots
\ea
converges to the unique fixed point $v$ of $T$ defined by $Tv=v$. 
\end{theorem}
We also paraphrase the ``local convergence'' result for Newton iteration.
\begin{theorem}[\cite{AHan2}, Thm.~.5.4.1, p236]
\label{th:newton}
Let $T:V\to V$ be Fr\'echet differentiable and let $T(U_*)=0$. Let also (*) $T'(U_*)^{-1}$ exist be a continuous linear map, and let $T'(U)$ be locally Lipschitz continuous in some neighborhood of $U_*$.  Then the Newton's iteration: given $U^{(o)}$, iterate
\ba
\label{eq:newton}
\Um=\Umm - T'(\Umm)^{-1}T(\Umm)
\ea
converges provided (**) $U^{(0)}$ is close enough to $U_*$. 
\end{theorem}

In many instances the functions we work with do not have continuous or Lipschitz derivatives. Thus we replace $T'$ in \eqref{eq:newton} by some selection from its subgradient $\partial T$. Weaker hypotheses in \cite{Ulbrich} allow to work in the {\em semismooth} framework. We paraphrase some results below as they are used for $T(U_*)=0$.
\begin{theorem}[\cite{Ulbrich}, Prop. 2.12, p29, Prop. 2.26, p35]
\label{th:snewton}
Let $\Hi=\R^M$, and let $T:\Hi \to \Hi$ be continuous and piecewise smooth. Let also (*) all subgradients from $\partial T(U_*)$ be nonsingular in some neighborhood of $U_*$ \mpc{and let $U^{(0)}$ be sufficiently close to $U_*$}.  Then the Newton's iteration \eqref{eq:newton} converges. 
\end{theorem}
%We will apply this result many times in what follows.

%%%%%%%
\subsection{Evolution equation on $\Hi$ and finite difference approximation}
We consider a nonlinear evolution equation on $\Hi$ with some $G:\Hi\to \Hi$
\ba
\label{eq:odem}
\ddt u(t) + G(u) =0, u(0)=\uz \in \Hi.
\ea
Let $\tau>0$ be the time step, and let $t_n=n \tau$.  The fully implicit finite difference approximation \mpc{of \eqref{eq:odem}} is 
\ba
\label{eq:fd}
\frac{1}{\tau}(U^n-U^{n-1}) + G(U^n)=0
\ea
where $\mpc{U^n} \approx u(t_n)$. The scheme \eqref{eq:fd} is of first order in $\tau$ if $G(\cdot)$ is smooth enough. 

It follows from \rthe{th:AHAN} that \eqref{eq:fd} rewritten as 
\ba
\label{eq:fdG}
U^n+ \tau G(U^n)=U^{n-1}
\ea
has a unique solution at each time step $t_n$ under some assumptions.

\begin{lemma}
\label{lem:trickG}
Let (i) $G$ be Lipschitz and monotone, or (ii) $G$ be Lipschitz and $\tau$ small enough. Then \eqref{eq:fdG} has a unique solution.
\end{lemma} 
\begin{proof}
For the case (i) the proof follows immediately  by \rthe{th:AHAN} since then $I+\tau G$ is strongly monotone. Consider (ii) when $G$ is not monotone. For sufficiently small $\tau$, $(I+\tau G)$ is strongly monotone: expanding $\inner{(I+\tau G)(u_1)-(I+\tau G)(u_2),u_1-u_2}$
\begin{multline}
\inner{u_1-u_2,u_1-u_2} + \tau \inner{G(u_1)-G(u_2),u_1-u_2} 
\\
= \norm{u_1-u_2}{}^2 - \tau \inner{G(u_1)-G(u_2),u_2-u_1}.
\end{multline}
The second term can be bounded from below by $-\tau L_G \norm{u_1-u_2}{}^2$ since 
$\inner{G(u_1)-G(u_2),u_1-u_2} \leq L_G \norm{u_1-u_2}{}^2$. With small $\tau$ so $c_T=1-\tau L_G>0$ and $(I+\tau G)$ is strongly monotone. 
\end{proof}

Formally we write the solution to \eqref{eq:fd} as $U^n=(I+\tau G)^{-1}U^{n-1}$, even if we do not form $(I+\tau G)^{-1}$ explicitly. In addition, if $G$ is nonlinear,  solving \eqref{eq:fd} requires an iteration.  For example, we can set up  fixed point iteration at every time step using previous time step 
\ba
\label{eq:fdfixed}
U^{n,(m)}+ \tau G(U^{n,(m-1)})=U^{n-1}, m=1,2,...; \;\; U^{n,(0)}=U^{n-1}.
\ea
By \rthe{th:fixed}, the iteration converges if the map $T$ in $ U^{n,(m)}=T(U^{n,(m-1)})=
U^{n-1}-\tau G(U^{n,(m-1)})$ is a contraction which follows if $\tau L_G$ is small enough.

We can also set-up Newton's iteration \eqref{eq:newton} for $T(U)=U+\tau G(U)-U^{n-1}$. Here, we might need to require that $T'(U)=I+\tau G'(U)$ exists and satisfies at least (*) from \rthe{th:newton}. Usually if $\tau$ is small enough, the initial guess from the previous time step solution $U^{n,(0)}=U^{n-1}$ suffices for  (**). 

The proofs in this paper will be similar to those we outline for \eqref{eq:fdG}.

%%%%
\subsection{Evolution equation under constraints on $\R$ and generalized play model of hysteresis}
\label{sec:evolution}
The results below are needed to develop the hysteresis models. We introduce the notion of constraint graphs on $\R$, and of evolution under constraints, well developed in \cite{Brezis73}. We refer also to \cite{Show-monotone}, and to our recent work on adsorption with hysteresis in \cite{P97,PS98,PS20,PS21} using linear play, nonlinear play, and generalized play models of hysteresis for the adsorption applications. In this paper we use generalized play models. 

%%%%%%%%%%%%%%%%%%%%%%%%%%%%%%%
\subsubsection{Constraint graph on $\R$}

Let $a \leq b \in \R$, and consider the non-empty closed interval $[a,b] \subset \R$. If $a < b$ we define the set-valued relationship
\ba
\label{eq:defr}
\cab{s} = \begin{cases}
(- \infty,0] \text{ if } s = a, \\
\{0\} \text{ if } a < s < b,  \\
[0,\infty) \text{ if } s = b.  \end{cases}
\ea
This definition indicates that $\cab{\cdot}$ is set-valued with domain $\Dom(\cC(a,b;\cdot))=[a,b]$; its graph is denoted by $\ccab=\{a\} \times (-\infty,0] \cup (a,b) \times \{0\}  \cup \{b\} \times [0,\infty]$. The relation $\cab{}$ is {\em monotone:} $\inner{v_1-v_2,w_1-w_2}_\Hi \geq 0$ for any $v_j,w_j\in \R$: $w_j \in \cg(v_j)$. It is also {\em maximal monotone} because the range $\Rg(\cg+I)=\R$,
and $\Rg(\lambda \cg+I)=\R$ for all $\lambda > 0$.  The resolvent 
$\rab{s} = (I+\tau \cab\cdot)^{-1}(s)$ is a function defined on $\R$ which can be written explicitly
\ba
\label{eq:rdef}
\rab{s} = (I+\tau \cab\cdot)^{-1} = \begin{cases}
a \text{ if } s \leq a, \\
s \text{ if } a <s < b, \\ 
b \text{ if } s \geq b,
\end{cases}
\; s \in \R,
\ea
and is independent of $\tau > 0$. The resolvent $\rab$ is used to determine the solution to a stationary problem
%%%%
\ba
\label{eq:stationary}
v + \tau \cab{v} \ni f \in \R,
\ea
in which the  symbol $\ni$ indicates that the left hand side of \eqref{eq:stationary} is a set. But once $v=\rab{f}$ is found, the selection $\cab{v} \ni  c^*=\tfrac{1}{\tau}(f-v)$ is unique.  

The following result is straightforward.
\begin{lemma}
\label{lem:resolvent}
For any $\alpha<\beta$ the resolvent \mpa{$\rab{s}$} 
is monotone in $s$ and Lipschitz continuous in $s$ with $L_R=1$, i.e, we have $\forall \xi_1,\xi_2 \in \R$
\bsa
\ba
\label{eq:RL}
\abs{\rab{\xi_1}-\rab{\xi_2}} \leq \abs{\xi_1-\xi_2)},
\\
0 \leq \inner{\xi_1-\xi_2,\rab{\xi_1}-\rab{\xi_2}} \leq \abs{\xi_1-\xi_2}^2. 
\ea
\esa
%%%
$R$ is also differentiable except at $a,b$ with $0 \leq R' \leq 1$ at the points of differentiability.
We also have from \mpcitee{Ulbrich}{Ex.2.1, p33} that $R$ is semismooth. 
\end{lemma}

%%%%%%%
\subsubsection{ODE with a constraint graph} 
Now we consider the Cauchy problem
\ba
\label{eq:ode}
\tfrac{d}{dt}v(t)+\cab{v(t)} \ni f(t),\ t \in (0,T],\  v(0) = \vz \in [a,b]. 
\ea
Its solutions are defined as limits of step functions $v_{\tau}$ built from implicit finite difference approximations $V^n \approx v(t_n)$ which solve
%%%%%%%%%%%%%%
%
\ba \label{eq:oden}
\tfrac{V^n - V^{n-1}}{\tau} + \cab{V^n} \ni f^n, \quad 1 \le n \le N,\ V^0=\vz,
\ea
for some suitable $f^n \approx f(t_n)$. Note that $V^n$ is well defined since \eqref{eq:oden} has the form  \eqref{eq:stationary} and thus we can use the resolvent
\ba
\label{eq:resolvent}
V^n=\rab{V^{n-1}+\tau f^n}.
\ea
\mpc{In turn, }the selection $c^n\in \cab{V^n}$ given \mpc{from \eqref{eq:oden} as } $\tau c^n=V^{n-1}+\tau f^n-V^n$ is the unique element of minimal norm \ourcite{Brezis73}{pp. 66, 28} \mpc{for which \eqref{eq:oden} holds}.

Based on the assumption that $\vz \in [a,b]$ and $f\in L^1(0,T)$, one can prove convergence of the finite difference solutions to $v(t)$; the rate is usually $O(\tau)$ in a Hilbert space. More details including the construction of $v_{\tau}$ and regularity of the solutions can be found in \cite{Brezis73,Show-monotone}. 

The construction \eqref{eq:resolvent} for the approximation to \eqref{eq:ode} is fundamental in the models in this paper. We illustrate this convergence in Section~\ref{sec:ODE}.

%%%
\subsubsection{ODE with a time-dependent constraint in generalized play hysteresis}
\label{sec:hysteresis}
Now we let the constraints in $\cab\cdot$ vary in time. To distinguish from fixed $a\leq b$, we consider the graph $C(\alpha,\beta;\cdot)$ with $\alpha(t) = \fr(t) \leq \beta(t) = \fl(t))$, with the resolvent $R(\alpha,\beta;\cdot)$. For more modeling flexibility, $\alpha,\beta$ can depend on some input $u(t)$,  namely, $\alpha(t) = \fr(u(t))$, $\beta(t) = \fl(u(t))$. Here $\fr,\fl$ are {continuous} monotone functions with $\fr\leq \fl$. 

Given input $u(t)$ and $\fr,\fl$, we generalize \eqref{eq:ode} to the Cauchy problem
\ba
\label{eq:odegh}
\tfrac{d}{dt}v+\cg(\fr(u(t)),\fl(u(t));v) \ni 0,\ t \in (0,T],\  v(0) = \vz. 
\ea
The approximation of \eqref{eq:odegh} follows similar to that for \eqref{eq:oden} but requires that we know $U^1,U^2,\ldots$. Also, we require $\vz \in [\fr(U^0),\fl(U^0)]$. We define the approximation $V^n$ as the solution to
\ba
\label{eq:genn}
\frac{1}{\tau}(V^n-V^{n-1})  + \cg(\fr(U^n),\fl(U^n);V^n) \ni 0\,,
\quad V^0 = \vz\,.
\ea
Given $U^n$, the solution $V^{n}$  follows using the resolvent \eqref{eq:resolvent}.
\ba
\label{eq:vsol}
V^n = \resab{\fr(U^n)}{\fl(U^n)}{\vbar} = \begin{cases}
\fr(U^n)\text{ if } \vbar \leq \fr(U^n), \\
\vbar \text{ if } \fr(U^n) <\vbar < \fl(U^n), \\ 
\fl(U^n) \text{ if } \vbar \geq \fl(U^n).  \end{cases} \vbar=V^{n-1}.
\ea
%%

%%%%%%%%%%%%%%%%%
\section{PDE model  with equilibrium, non-equilibrium, and hysteresis variants and its fully discrete approximation}
\label{sec:model} 
In this section we present details of the model \eqref{eq:uchi}. We closely follow our presentation in \cite{BPV} and recent in \cite{VP-Tp,PHV} which we augment with the non-equilibrium or hysteretic close for $\chi(u)$.  Our equilibrium model is similar to those presented in the applications literature; e.g., in 
\cite{RomanovskyOsterkamp2000,Michalowski93,NicolskyRomanovskyTipenko07}, but we make simplifying assumptions and do not consider any particular physical scenarios.

We start with the presentation of an equilibrium model in physical quantities in Section~\ref{sec:modelph}. Next we reformulate \mpc{that model} in simplified notation in Section~\ref{sec:modelu}. We follow up with the non-equilibrium models, and their discretization. 

%%%
\subsection{Physical model}
\label{sec:modelph}
We consider the energy conservation in a partially frozen soil at the Darcy scale, where the soil is treated as a continuum mixture of solid rock grains and water in the liquid or ice phase. The void space available to the water has volume fraction $\eta$ also called the porosity, and the total energy density is given by the enthalpy variable $W$; its evolution is balanced by the heat flux $\QT$ and the sources. The model is closed by Fourier's heat conduction \mpc{with a} nonlinear heat conductivity $\LambdaT$ \mpc{and reads}
\bsa
\label{eq:physicaltheta}
\ba
\label{eq:TD}
\partial_t ( W(\theta))
+\nabla \cdot \QT=0, \QT=-\lambdaT \nabla \theta.
\ea
The definitions of $W(\theta),\LambdaT(\theta)$ we adopt in this paper rely on the following modeling assumptions.  

{\em Physical assumptions:} The porous medium is rigid and homogeneous, i.e. the porosity $\eta$ is constant, and the only fluid present is water in the ice and liquid phases. The following physical material properties of the liquid, ice, and rock grains are constant including the densities $\rho_l,\rho_i,\rho_r$, the heat\ capacities\ $ c_l,c_i,c_r$, and the heat\ conductivities $k_l,k_i,k_r$, respectively. We also assume the latent heat $L$ is constant. 

Denoting by $S(\theta)$ the volume fraction of liquid, $W$ includes the volume fraction weighted thermal energy densities plus the latent heat for the water material. Given $\eta$, and denoting $c_u = \eta c_l + (1-\eta) c_r$ and $c_f = \eta c_i + (1-\eta) c_r$ we have as in \cite{BPV}
\ba
\label{eq:Alpha}
W&=& \int_{0}^{\theta} (c_u S(v) + c_f (1-S(v))dv + L \eta S.
\ea
We also require $\LambdaT$ which can be found by upscaling as in \cite{PVB} or with simple parametric volume fraction weighting. Wlog, in this paper we assume arithmetic weighting. \mpc{We denote} $k_u = \eta k_l + (1-\eta) k_r$ and $k_f = \eta k_i + (1-\eta) k_r$, \mpc{and define}
\ba
\label{eq:LambdaT}
\LambdaT=
%(1-\eta)k_r+\eta(k_lS+k_i(1-S)) = 
k_uS+k_f(1-S)=(k_u-k_f)S+k_f.
\ea

To complete the model we require $S(\theta)$. In equilibrium it is described  by an algebraic expression \mpc{which} we denote by $F(\cdot)$. \mpc{In this paper we} adopt  a simplification of the parametric model \cite{Michalowski93} 
\ba
\label{eq:chiM}
S=F(\theta) = \left\{\begin{array}{ll}
1; & \theta \geq 0,\\
e^{\mpc{b \theta} }; & \theta < 0.
\end{array}\right.
\ea
\esa
Here $b>0$ can be found empirically or from upscaling as in \cite{PHV}. The physical model \eqref{eq:physicaltheta} is now complete and only requires initial and boundary conditions. 

%%%%%%%%%%%%
\subsubsection{\bf Simplified PDE model} 
\label{sec:modelu}
We derive now a simplified model which retains the qualitative character of \eqref{eq:TD} but uses rescaled primary unknowns and data, \mpc{under the assumption that $\eta$ is constant}.  To distinguish between the physical set-up and the simplified model, we replace $\theta$ by $u$. We set $\tc_u=\tfrac{c_u}{\eta L},  \tc_f=\tfrac{c_f}{\eta L}, \tk_u=\tfrac{k_u}{\eta L}; \tk_f=\tfrac{k_f}{\eta L}$, and divide \eqref{eq:TD} by $L\eta$. We have
\bsa
\label{eq:physical}
\ba
\tfrac{W(u)}{\eta L} &=& c(u)+S = (\tc_u-\tc_f)\int_o^u S(v)dv+\tc_fu +S.
\\
\tfrac{\LambdaT}{L \eta}&=&k(u)  = \tk_f+(\tk_u-\tk_f) S.
\ea
\esa
With these, if we replace $S$ by $\chi$, the model reads the same as \eqref{eq:uchi}
\bsa
\label{eq:pdeall}
\ba
\label{eq:pde}
\partial_t ( c(u)+\chi)
-\nabla \cdot(k(u)\nabla u)=\tilde{f}.
\\
u(x,0)=\uz(x),
\\
u(x,t)\vert_{\partial \Omega}=u_D(t),
\ea
\esa
where $u_D(t),\uz(\cdot)$ as well as the rescaled source term $\tilde{f}$ are given. 
This model will be next supplemented by (EQ) or (NEQ) or (HYST) relationships which define the relationship $(u,\chi)$; in general, $\chi$ is not the same as $F(u)$ and has its own dynamics. 

One must make a choice how $k(u),c(u)$ defined in \eqref{eq:physical} (depending on $S$) is evaluated. One option is to use the current value of $S=\chi$. Another simpler option is to substitute the equilibrium relationship $S=F(u)$ in \eqref{eq:physical}.  

\begin{remark}
In this paper we make the following choice. To calculate $c(u), k(u)$, we plug in the equilibrium formula $S=F(u)$ \eqref{eq:chiM} in \eqref{eq:Alpha} and \eqref{eq:LambdaT} rather than recalculate $c(u)$ and $k(u)$ based on the current $\chi$, an unknown in the (NONEQ) and (HYST) models. This approach allows to focus on the abundant remaining challenges. 
\end{remark}

\begin{example}
\label{ex:ck}
We take $S=F(u)$ given in \eqref{eq:chiM} and  and evaluate $c(u),k(u)$ 
\bsa
\label{eq:kc}
\ba
c(u) =(\tc_u-\tc_f)\int_o^u F(v)dv+\tc_fu = \begin{cases}\frac{\tc_u-\tc_f}{b} (e^{bu}-1)+\tc_f u,&u\leq 0,\\\tc_u u,&u>0,\end{cases}
\\
k(u)=\tk_f+(\tk_u-\tk_f) F(u)=\begin{cases}(\tk_u-\tk_f) e^{bu}+\tk_f,&u\leq 0,\\k_u,&u>0. \end{cases}.
\ea
\esa
\end{example}
\begin{remark}
\label{rem:cimplicit}
We can comment now on the behavior of $c(u),k(u)$ in \eqref{eq:kc}. Both are continuous but neither is differentiable at $u=0$. Therefore it is not recommended in numerical models to replace $\partial_t c(u)$ by $\frac{dc}{du}\partial_t u$ since $c'(u)$ features a large  jump and causes numerical oscillations. Implicit treatment of $c(u)$ is recommended.
\end{remark}

\begin{remark}
We can absorb the constants $L\eta$ in \eqref{eq:physical} making the time variable nondimensional to adhere to the characteristic time of thawing/freezing of about $\mpunit{day}=8.64 \times 10^4\mpunit{s}$. See our simulations in Section~\ref{sec:results-pde}.
\end{remark}

%%%%%%%%%%%%%%%%%%%%
\begin{assum}
\label{ass:ckF}
We will make the following assumptions on the data $c,k$ and $F$. 
\\
{\bf (AC):} The energy density $c(u)$ is continuous monotone increasing in $u$ and smooth except at $u=0$ \mpc{where} $c(0)=0$. It is also differentiable except at $u=0$, with $L_c \geq c'(u)\geq c_{min}>0$.
\\
{\bf (AK):} Heat conductivity $k(u)$ is bounded above and below by positive constants, is smooth except at $u=0$, and Lipschitz continuous with constant $L_k$.  
\\
{\bf (AF):} The equilibrium relationship $S=F(u)$ is nonnegative and bounded above, monotone nondecreasing in $u$, smooth except at $u=0$, and globally Lipschitz, i.e., there are nonnegative constants $L_F,F_{\infty},F_0$ as follows: for every $u_1,u_2 \in \R, u\in \R$ 
\bsa
\label{eq:Fprop}
\ba
\label{eq:Fmon}
(F(u_1)-F(u_2))(u_1-u_2) &\geq& 0,
\\
\label{eq:FLip}
\abs{F(u_1)-F(u_2)} &\leq& L_F \abs{u_1-u_2},
\\
\label{eq:Finf}
0 \leq F(u) &\leq& F_{\infty},
\\
F(0)=F_0&\leq& F_{\infty}.
\ea
\esa
\end{assum} 
The data $c(u),k(u)$ from Example~\ref{ex:ck} and $F(u)$ from \eqref{eq:chiM} satisfy these assumptions, and we have $L_F=b,F_0=F_{\infty}=1$. 

%%%
\subsection{Complete simplified model and fully discrete counterpart}

Now we write the spatially discrete approximations to the IBVP system \eqref{eq:pdeall}. For simplicity, in this paper we consider only the cell-centered finite differences, with some rectangular grid $\tT_h$ over $\Omega$, with cell centers $x_j,j=1,2 \ldots M$. 

The discrete system is solved for the approximations $U_j^n \approx u(x_j,t_n), \Chi_j^n \approx \chi(x_j,t_n)$ with the source terms denoted by $f_j^n$. The corresponding time continuous variables are, respectively, $U_j(t) \approx u(x_j,t)$ and $\Chi_j(t) \approx \chi(x_j,t)$. 
These are collected in the vectors $U^n,\Chi^n$, and $U(t),\Chi(t)$, respectively.  The choice of these approximations on a uniform grid does not require the use of non-diagonal mass matrices, \mpc{thus it makes the notation easier  and allows  }to focus on the features of (non)equilibria. Also, with the CCFD approximation, the boundary data $u_D(t)$ is absorbed in the right hand side $f^n_j$, and there are no degrees of freedom on the boundary $\partial \Omega$ which is assumed to align with the edges of elements in $\tT_h$. We refer to \cite{PJW02,BPV} for extensive details on the discretization, and to \cite{APS} for recent analysis involving nonlinear diffusion terms similar to $-\nabla \cdot (k(u) \nabla u)$. 

\begin{remark}
The node-centered FD approximation can be handled with the same analysis we employ below, as long as it is done for the interior unknowns only.  
\end{remark}

The continuous in time discrete in space approximation to \eqref{eq:pdeall} is an ODE in $\Hi=\R^M$ and reads
\ba
\label{eq:pdesystC}
\ddt( C(U)+\Chi)
+A(U)U=f(t); \; (C(U)+\Chi)\vert_{t=0}=W_{init}.
\ea

Here the diffusion matrix $A=A(U)$ because $k=k(u)$ is nonlinear, and the entries of $A$ are linear in $k$.  Thus the relationship \mpc{$U \to A(U)$} is Lipschitz because of Assumption (AK).  {However, $L_A$ depends on $L_{k}$ and the dimension $M$.}

\begin{remark}
Based on Assumption~\ref{ass:ckF} 
we see that $\inner{C(U)-C(V),U-V}\geq c_{min}\norm{U-V}{}^2$. Also, since $c(0)=0$, we have $\inner{C(U),U}\geq c_{min}\norm{U}{}^2$, i.e. $C(U)$ is strongly monotone. It is also Lipschitz. Thus it plays a similar role to $I$ in analysis. We shall  thus   consider $C(U)=U$ in what follows which will considerably simplify the notation.  
\end{remark}

After this simplification we consider
\ba
\label{eq:pdesyst}
\ddt( U+\Chi)
+A(U)U=f(t); \; (U+\Chi)\vert_{t=0}=W_{init}.
\ea
The fully discrete version of \eqref{eq:pdesyst} is
\ba
\label{eq:pdesysn}
\tfrac{1}{\tau}( U^n+\Chi^n)
+A(U^n) U^n=f^n; (U^0+\Chi^0)=W_{init}.
\ea
Next we must add to \eqref{eq:pdesysn} some relationship binding $U$ and $\Chi$ pointwise; specifically we add the discrete counterparts of \eqref{eq:EQ} \mpc{or} \eqref{eq:NEQ} \mpc{or} \eqref{eq:HYST}. This is done below.

%%%%%%%%%%%%%%%%%%
\subsubsection{Equilibrium system based on \eqref{eq:EQ}}
\label{sec:EQn}
In equilibrium, we have \eqref{eq:pdesysn} and discrete form of \eqref{eq:EQ}
\bsa
\label{eq:eqsysn}
\ba
\label{eq:eqsysnu}
\tfrac{1}{\tau}( U^n+\Chi^n)
+A(U^n) U^n &=& f^n;
\\
\label{eq:eqsysnv}
\Chi^n=
(\Chi_j^n)_j&=&(F(U_j^n))_j
=\fF(U^n).
\ea
\esa
With the pointwise relationship \eqref{eq:eqsysnv} \mpc{at every degree of freedom, we see that $U \to \fF(U)$} is Lipschitz with $\LfF=L_F$ since $F$ is Lipschitz component-wise \eqref{eq:FLip}.

The system \eqref{eq:eqsysn} is nonlinear with two nonlinearities $\fF$ and $A$, and must be solved by iteration. Its analysis will be given in Section~\ref{sec:EQa}.

%%%%%%%%%%%%%%%%%%%
\subsubsection{Non-equilibrium model based on \eqref{eq:NEQ}}
\label{sec:NEQn} 
Now we discuss the non-equilibrium relationship \eqref{eq:NEQ}. Such relationships are fairly common in systems with large spatial scales. In principle, the model \eqref{eq:uchi} is postulated at any spatial scale $x$ and time scale $t$. Then the temperature  $u(x,t)$ as an intensive variable  is well understood as the pointwise quantity \mpc{ and the} formation and disappearance of ice phase is in complete equilibrium with $u$. This understanding \mpc{requires that} we can resolve the fine details of the phase behavior, i.e, we work at a microscopic scale at which the ice/liquid interfaces are visible. 

However, it is in general difficult to approximate the solutions to \eqref{eq:uchi} \mpc{at that scale}. \mpc{Therefore $u(x,t)$ should}  
be understood as an average \mpc{at time $t$} of the temperature over some local region $\omega(x)$. Then it is natural that in any such region one can observe simultaneously both the ice crystals and liquid regions. 
\begin{remark}
\label{rem:cup}
A simple example of non-equilibrium is provided by a cup of water to which we add a few cubes of ice. If $u(t)$ is the average temperature in the entire cup at time $t$, then one can 
talk about an equilibrium between $u$ and the average water fraction $S=\chi(u)$ in the cup only after sufficiently large time. The process of getting to the equilibrium is modeled by \eqref{eq:NEQ} with rate $B>0$. 
\end{remark}

We rewrite now \eqref{eq:pdesysn} and  finite difference approximation to \eqref{eq:EQ}. The fully implicit finite difference scheme reads, after some rearranging and setting $\bar{B}=\frac{1}{1+\tau B}$.
\bsa
\label{eq:neqsysn}
\ba
\label{eq:neqsysnu}
\frac{1}{\tau}(U^n-U^{n-1})  + \frac{1}{\tau}(\Chi^n-\Chi^{n-1}) +A(U^n)U^n= f^n,
\\
\label{eq:neqsysnv}
\Chi^n=\frac{B\tau}{1+\tau B} \fF(U^n)+ 
\frac{1}{1+\tau B}\Chi^{n-1} 
%\\
%\nonumber
= \fF(U^n) (1-\bar{B})+ 
\bar{B}\Chi^{n-1}.
\ea
\esa

Comparing \eqref{eq:neqsysnv} to \eqref{eq:eqsysnv} we see that if $\tau B>>1$, then $\bar{B}$ is small. \mpc{Also,}  
\eqref{eq:neqsysnv} and \eqref{eq:eqsysnv} have similar right hand sides, with a small role played by $V^{n-1}$.  

From this observation we see that the solvability of the system \eqref{eq:neqsysn} follows similarly to that of \eqref{eq:eqsysn} once we plug \eqref{eq:neqsysnv} to \mpa{\eqref{eq:neqsysnu}}. Details will be given in Section~\ref{sec:NEQa}. 

%%%%%%%
\subsubsection{Hysteresis model with \eqref{eq:HYST}}
\label{sec:HYSTn}
Now we generalize the model in Section~\ref{sec:NEQn} to acknowledge the fact that the freezing processes have characteristics distinct from thawing. In particular, continuing Remark~\ref{rem:cup}, supplying the heat to the cup will eventually result in all water in the liquid state with rate $B_{thaw}$. However, if the cup is placed in a freezer, the water will gradually turn to ice, but this process due to the slow nucleation of ice crystals  will likely proceed with a different rate $B_{freeze} \neq B_{thaw}$, and it might even follow a different $F_{freeze}(u)$ as $u$ is decreasing than $F_{thaw}$ as $u$ is increasing. This is the simplest conceptual model of hysteresis phenomena from physical standpoint.  

Hysteresis models have been well studied \cite{Visintin94,KrasPokr89,Mayergoyz}. In this paper we use the concept of generalized play models of hysteresis recently discussed in the practical setting in \cite{PS21} and given in Section~\ref{sec:hysteresis}. We now apply it to the specific characteristics of hysteresis in the permafrost. 

We wish to use a particular generalized play model \eqref{eq:odegh} relating $\chi$ to $u$. We must now
make a particular selection of $\fr(u),\fl(u)$ to fit the model to a permafrost application. We shall consider a given equilibrium function $F$ satisfying (AF), and a function $G$ satisfying (AF) which together satisfy $G \geq F$, and additional conditions that will be relevant later. In particular, we wish to model 
\ba
F(u) \leq \chi \leq G(u),
\ea
which we rewrite $0 \leq \chi-F(u) \leq G(u)-F(u)$. In other words, we set up $\alpha=\fr=0$ and $\beta=\fl(t)=G(u(t))-F(u(t))$, and apply \eqref{eq:odegh} 
\ba
\label{eq:hystFG}
\ddt \chi + C(0,\beta;\chi-F(u))\ni 0; \chi(x,0)=\chi_{init}(x).
\ea
This equation now supplements \eqref{eq:pdeall}.

With this specific choice we apply the scheme  \eqref{eq:genn} with some discrete value $\beta_j^n$ known at every degree of freedom $j$, to obtain 
\ba
\label{eq:hystFGn}
\frac{1}{\tau}(\Chi_j^n-\Chi_j^{n-1})  + \mpc{\cg(0,\beta_j^n;\Chi_j^n-F(U_j^n))} \ni 0\,; \; \Chi_j^0=\Chi_{init,j}.
\ea
Now we rearrange, subtracting $F(U_j^n)$ from both sides
\ba
\label{eq:chibeta}
\Chi_j^n-F(U_j^n)  + \tau \cg(0,\beta_j^n;\Chi_j^n-F(U_j^n)) \ni \Chi_j^{n-1}-F(U_j^n)\,.
%\\
%\Chi^n-F(U^n)  + \tau \cg(0,G(U^n)-F(U^n);\Chi^n-F(U^n)) \ni \Chi^{n-1}-F(U^n)\,,
\ea
Applying the resolvent \eqref{eq:resolvent} we get 
$\Chi_j^n-F(U_j^n) =(R(0,\beta_j^n;\vbar_j))_j$ with $\vbar_j=\Chi_j^{n-1}-F(U_j^n)$ 
\ba
\label{eq:vsolF}
\Chi^n = \fF(U^n)+(\resab{0}{\beta_j^n}{\vbar_j})_j; \Chi^n_j=
%%%
\begin{cases}
F(U_j^n)\text{ if } \vbar_j \leq 0, \\
F(U_j^n)+\vbar_j \text{ if } 0 <\vbar_j < \beta_j^n, \\ 
\beta_j^n \text{ if } \vbar_j \geq \beta_j^n.  \end{cases}
\ea
It remains to specify $\beta^n$. 

The fully implicit choice \mpc{where} pointwise  $\beta_j^n=G(U_j^n)-F(U_j^n)$ seems natural, but requires iteration and is challenging in analysis. 

In this paper we consider therefore the time lagging (sequential approach) with $\beta_j^n=G(U_j^{n-1})-F(U_j^{n-1})$. 
The fully implicit finite difference scheme for the system \eqref{eq:pdeall}, \eqref{eq:hystFG} reads now
\bsa
\label{eq:hystsysn}
\ba
\label{eq:hystsysnu}
\frac{1}{\tau}(U^n-U^{n-1})  + \frac{1}{\tau}(\Chi^n-\Chi^{n-1}) +A(U^n)U^n\ni (f_j^n)_j,
\\
\label{eq:hystsysnv}
\Chi^n = \fF(U^n) + (\resab{0}{\beta_j^n}{\Chi_j^{n-1}-F(U_j^n)})_j,
\\
\nonumber 
\beta_j^n=G(U_j^{n-1})-F(U_j^{n-1}).
\ea
\esa

The analysis of this system is given in Section~\ref{sec:HYSTa}. We see it will have some similarity with that of \eqref{eq:neqsysnv} 
since the right hand side in \eqref{eq:hystsysnv} is somewhat similar to that in 
\eqref{eq:neqsysnv} and \eqref{eq:eqsysnv}, with a small role played by $V^{n-1}$. 

%%%%%%%%%%
\section{Well-posedness of discrete systems and convergence of iterations}
\label{sec:analysis}
Now we study systems \eqref{eq:eqsysn}, \eqref{eq:neqsysn}, \eqref{eq:hystsysn} as variants of 
\ba
\label{eq:nu}
U^n+ \Chi^n +
\tau A(U^n) U^n &=& g^n =\tau f^n + U^{n-1}+\Chi^{n-1},
\ea
coupled with a relationship for $\Chi^n$ in terms of $\Chi^{n-1}$ and $\fF(U^n)$ specific to the particular model (EQ), (NEQ), or (HYST). 

Because of the double nonlinearity in $\fF$ and $A$, we start with introductory material on handling $A(U)U$ to be followed by the material on $\fF(U)$.  

Next we study well-posedness as well as a fixed point iteration; we give proofs for the (EQ) model \eqref{eq:eqsysn}, (NEQ) model \eqref{eq:neqsysn}, and (HYST) model \eqref{eq:hystsysn} in Sections~\ref{sec:EQa}, \ref{sec:NEQa}, and \ref{sec:HYSTa}, respectively. 

\subsection{Case of trivial $\Chi^n=0$ in \eqref{eq:nu}} 
We first set $\Chi^n=0$ in \eqref{eq:nu} and solve
\ba
\label{eq:un}
U^n+
\tau A(U^n) U^n &=& g^n.
\ea
\begin{theorem}
\label{th:u}
Let $\Hi=\R^M$. For any $u \in \Hi$, let $A(u)$ be a symmetric uniformly positive definite and uniformly bounded operator which is Lipschitz in $u$
\bsa
\label{eq:Aprop}
\ba
\label{eq:Aspd}
\inner{A(u)\xi,\xi} \geq \kappa_0 \norm{\xi}{}^2,
\\
\label{eq:Abd}
\inner{A(u)\xi,\zeta} \leq A^{max}\norm{\xi}{}\norm{\zeta}.
\\
\label{eq:ALip}
\norm{(A(u)-A(v))\xi} \leq L_A \norm{u-v}{}\norm{\xi}{}
\ea
\esa
Then if $\tau$ is small enough, there exists a unique solution to  
\ba
\label{eq:u}
U+\tau A(U) U &=& g,
\ea
which satisfies the bound 
\ba
\label{eq:ubound}
\norm{U}{} \leq 
\tfrac{1}{ (1+\tau \kappa_0)} \norm{g}{}.
\ea
\end{theorem}
%%%
\begin{proof}
We set up fixed point iteration. Given some initial guess $U^{(0)} \in \Hi$ we iterate
\ba
\label{eq:m}
U^{(m)}+\tau A(U^{(m-1)}) U^{(m)} &=& g, m=1,2,\ldots
\ea
Fix $m$ now. Taking inner product with $U^{(m)}$ of \eqref{eq:m} and applying \eqref{eq:Aspd} we get
\ba
\label{eq:uest}
(1+\tau \kappa_0) \norm{U^{(m)}}{}^2 \leq 
\norm{U^{(m)}}{}^2 
+\tau \inner{A(U^{(m-1)})U^{(m)},U^{(m)}}=
\inner{g,U^{(m)}},
\ea
with the right hand side bounded by $\norm{g}{}
\norm{U^{(m)}}{}$. Dividing both sides by $\norm{U^{(m)}}{}$ we obtain the bound
\ba
\label{eq:umbound}
\norm{U^{(m)}}{} \leq \tfrac{1}{ (1+\tau \kappa_0)} \norm{g}{}, \; \forall m=1,2,\ldots.
\ea

Next we study the map $U^{(m-1)} \to U^{(m)}=T(U^{(m-1)})$ given by \eqref{eq:m}, and show $T$ is a contraction. To this aim we subtract \eqref{eq:m} for $m$ and that for $m-1$ 
%%%
\bas
\Umm+\tau A(\Ummm) \Umm &=& g,
\eas
and take the inner product of the result with 
$U^{(m)}-U^{(m-1)}$. We get, after adding and subtracting  $A(\Umm)\Umm$ inside the inner product
\begin{multline}
\label{eq:mAest}
\inner{\Um-\Umm,\Um-\Umm}
\\
+ \tau \inner{A(\Umm)\Um-A(\Ummm)\Umm,\Um-\Umm}]
\\
=\inner{\Um-\Umm,\Um-\Umm}
+ \tau \inner{A(\Umm)(\Um-\Umm),\Um-\Umm}
\\
+ \tau \inner{A(\Umm)\Umm - A(\Ummm)\Umm,\Um-\Umm}=0.
\end{multline}
Now shifting the last term on the left hand side to the right and estimating from \eqref{eq:ALip} using Cauchy-Schwarz inequality we obtain 
\begin{multline*}
\tau \inner{(A(\Umm)- A(\Ummm))\Umm,\Um-\Umm} 
\\
\leq 
\tau L_A \norm{\Umm-\Ummm}{} \norm{\Umm}{} \norm{\Um-\Umm}{}.s
\end{multline*}
For the second term we apply \eqref{eq:Aspd} to get
\begin{multline*}
%\label{eq:mestimate}
(1+\tau \kappa_0) \norm{\Um-\Umm}{}^2
\\
\leq \tau L_A \norm{\Umm-\Ummm}{} \norm{\Umm}{} \norm{\Um-\Umm}{},
\end{multline*}
and dividing by $\norm{\Um-\Umm}{}$ we see
\bas
\norm{\Um-\Umm}{}
\leq \frac{\tau L_A}{(1+\tau \kappa_0) } \norm{\Umm-\Ummm}{} \norm{\Umm}{}.
\eas
Combining with \eqref{eq:umbound} for $\norm{\Umm}{}$ we finally obtain that the map $T$ has Lipschitz constant 
\ba
\label{eq:uLT}
L_T = \frac{\tau L_A}{(1+\tau \kappa_0) } \norm{g}{}.
\ea
From this we conclude that if $\tau$ is small enough, then $L_T<1$, and the iteration \eqref{eq:m} converges to 
the unique solution of \eqref{eq:u} which satisfies the desired bound \eqref{eq:ubound}. 
\end{proof}
\begin{corollary}
Applying \rthe{th:u} we see that at every time step the solution $U^{n}$ to \eqref{eq:un} exists and is unique. Moreover, the fixed point iteration \eqref{eq:m} converges if $\tau$ is small enough. A natural choice of the initial guess $U^{n,(0)}=U^{(n-1)}$. 
\end{corollary}
%%%
\begin{remark}
The technique we applied does not work in the \mpa{infinite} dimensional setting since there, in particular, one cannot assume that $A(U)$ is uniformly bounded or Lipschitz. 
\end{remark}

%%%%%%%%%%%
\subsection{
Well-posedness and convergence for the equilibrium problem \eqref{eq:eqsysn}}
\label{sec:EQa}
Now we consider \eqref{eq:nu} in which we set $\Chi^n=\fF(U^n)$. \mpc{We have}
\ba
\label{eq:equn}
U^n+ \fF(U^n)+
\tau A(U^n) U^n &=& g^n;
\ea
Now we fix and suppress $n$ and study  
\ba
\label{eq:equAF}
U+ \fF(U)+
\tau A(U) U &=& g;
\ea
This problem features double nonlinearity $\fF(U)$ and $A(U)U$. The challenge is that we do not have lower bounds on $\inner{\fF(U),U}$, only \eqref{eq:Fmon}, thus we have to shift the terms involving $\fF(\cdot)$ to the right hand side which results in some restrictive estimates involving $L_F$.  We give analysis, and consider alternatives.

\subsubsection{Handling only nonlinearity in $F$ while lagging nonlinearity in $A$}
\label{sec:mA}
We can apply \rthe{th:AHAN} to the problem 
\ba
\label{eq:equbar}
U+ \fF(U)+
\tau \bar{A} U &=& g;
\ea
in which $\bar{A}$ is fixed, and satisfies all the desired properties  \eqref{eq:Aprop} with a trivial $L_A=0$. Since $F(\cdot)$ is monotone and Lipschitz, and since $\bar{A}$ linear spd, we obtain that the operator $T(U)=U+\fF(U)+\tau \bar{A}U$ is  strongly monotone with constant $c_0={1+\tau \kappa_0}$ and Lipschitz with the (large) constant $L_T=1+\LfF+\tau \norm{\bar{A}}{}$. 
We obtain the following.

\begin{lemma}
\label{lem:mF}
Let $\bar{A}$ be constant matrix which satisfies \eqref{eq:Aspd}, \eqref{eq:Abd}, and let $F(\cdot)$ satisfy \eqref{eq:Fprop}. Then the problem 
\eqref{eq:equbar} is well-posed and 
\ba
\label{eq:uboundA}
\norm{U}{}\leq \frac{1}{1+\tau \kappa_0}\norm{g}{}.
\ea
\end{lemma}
Next we check if the iteration $\Umm \to \Um=T(\Umm)$ defined by
\ba
\label{eq:mF}
\Um+ \fF(\Umm)+
\tau \bar{A} \Um &=& g,
\ea
converges. To do so, we use estimates similar to those in the Proof of \rthe{th:u}. 
\begin{multline}
(1+ \tau \kappa_0) \norm{\Um -\Umm}{}^2  \leq 
\norm{\fF(\Umm)-\fF(\Ummm)}{} \norm{\Um - \Umm}{}
\\
\leq 
L_F \norm{\Umm - \Ummm}{} \norm{\Um - \Umm}{}.
\end{multline}
The Lipschitz estimate employed above works since 
\bas
\norm{\fF(U)-\fF(V)}{}^2 =\sum_j \abs{F(U_j)-F(V_j)}^2 \leq L_F^2 \sum_j \abs{U_j-V_j}^2.
\eas
Still, we require $\frac{L_F}{1+\tau \kappa_0}<1$ for convergence, and even if this holds, convergence may be slow. Thus Newton solver may be a better option for this particular iteration. 

To ensure \eqref{eq:newton}  can be applied, we must check the conditions in \rthe{th:newton} or \rthe{th:snewton}. Consider the Jacobian $\jJ=T'(U)=I+\fF'(U)+\tau \bar{A}$ for $T(U)=U+\fF(U) + \tau \bar{A} U-g$. Here $\fF'(U)$ is a diagonal matrix with the $j'th$ entry equal $F'(U_j)$. Now $\jJ$ is nonsingular for any $U$, since as we said above $U+\fF(U)+\bar{A}U$ is strongly monotone. However, since $F'$ is not continuous, $T'$ is only given piecewise, and we are not able to verify that $T'$ is Lipschitz or that $(T')^{-1}$ is continuous (see already the scalar case $M=1$).  However,  $\fF$ is semismooth, thus $\jJ$ is as well. Thus Newton iteration converges.

\begin{lemma}
\label{lem:snewton}
(Semismooth) Newton iteration converges for \eqref{eq:equbar}.
\end{lemma}

An alternative proof of well-posedness exploiting the fact that the monotone operator $I+\fF'+\tau \bar{A}$ is a subgradient is given in \cite{BPV}. 
\subsubsection{Handling the double nonlinearity $F$ and $A$ directly}  
\label{sec:FAdirect}

Now we wish to handle \eqref{eq:equAF} rather than \eqref{eq:equbar}. For this we derive additional estimates but they require rather unrealistic restrictive assumptions. We illustrate the difficulties. 
We set up the iteration combining \eqref{eq:m} and \eqref{eq:mF} to get $\Umm \to \Um=T(\Umm)$ defined by
\ba
\label{eq:mAF}
\Um+ \fF(\Umm)+
\tau A(\Umm)\Um &=& g.
\ea
We follow the proof of \rthe{th:u} adding the elements of the proof of Lemma~\ref{lem:mF}. 

First we obtain an upper bound for $\norm{\Um}{}$ by following calculations that led to \eqref{eq:uest} modifying its right hand side from $g$ to $g-\fF(\Umm)$ and estimating its norm by $\norm{g}{}+F_{\infty}$, from $
\inner{\fF(\Umm),\Um}\leq
\norm{\fF}{\infty}\norm{\Um}{}\leq F_{\infty} \norm{\Um}{}$  to obtain 
\ba
\label{eq:uestAF}
\norm{\Um}{}
\leq \frac{1}{1+\tau \kappa_0} (\norm{g}{}+F_{\infty}).
\ea
Handling next the question of convergence of the iteration and proceeding as in \eqref{eq:mAest} we have an additional term on its right hand side
$\inner{F(\Umm)-F(\Ummm),\Um-\Umm}$, similar as in \eqref{eq:mF}. Thus 
\begin{multline}
(1+\tau \kappa_0) \norm{\Um-\Umm}{}^2
\\
\leq \tau L_A \norm{\Umm-\Ummm}{} \norm{\Umm}{} \norm{\Um-\Umm}{} 
\\
+
L_F \norm{\Umm - \Ummm}{} \norm{\Um - \Umm}{}.
\end{multline}
Dividing by $\norm{\Um-\Umm}{}$ and incorporating \eqref{eq:uestAF}
\bas
\norm{\Um-\Umm}{}
\leq \frac{\tau L_A+L_F}{(1+\tau \kappa_0) }  (\norm{g}{} +F_{\infty}) \norm{\Umm-\Ummm}{},
\eas
we see that $T$ is a contraction under a rather restrictive condition that 
\ba
\label{eq:LTAF}
L_T=\frac{\tau L_A+L_F}{(1+\tau \kappa_0) }  (\norm{g}{} +F_{\infty}) <1,
\ea
and we conclude with the result as follows.
\begin{prop}
Let $F,A$ be as in \rthe{th:u} and \rlem{lem:mF}, respectively, and let \eqref{eq:LTAF} hold. Then a unique solution to \eqref{eq:equAF} exists and the iteration \eqref{eq:mAF} converges. 
\end{prop}

In general the assumption \eqref{eq:LTAF} cannot be verified. Therefore we proceed differently in the next section. 

\subsubsection{Handling double nonlinearity by double iteration}
\label{sec:double}
As mentioned above, solving \eqref{eq:equAF}  with one fixed iteration is not a practical option since \eqref{eq:LTAF} can be restrictive. On the other hand, since $U+\fF(U)$ is strongly monotone, there are better algorithms than the single monolithic iteration \eqref{eq:mAF}. 

We solve \eqref{eq:equAF} by double iteration exploiting \eqref{eq:equbar}. Given $\Umm$, we find $\Um$ as follows:
\ba
\label{eq:double}
\text{\ Calculate\ }\bar{A} = A(\Umm). 
\text{\ Solve\ } \Um+ \fF(\Um)+
\tau \bar{A} \Um = g;
\ea
using Newton's method. Here we iterate 
$U^{(m,0)}, U^{(m,1)}, \ldots U^{(m,r)}, \ldots$. 
In each iteration $r$ we have the Jacobian
\bas
\jJ=I+\fF'(U^{(m,r)})+\tau \bar{A}
\eas
%%%
which  is spd for any $U^{(m,r)}$. 
From Lemma~\ref{lem:snewton} we have therefore convergence of Newton's step within the double iteration. With a good initial guess from previous time step $U^{n,(0)}=U^{n-1}$, the iteration is likely very robust for the equilibrium problems, as we have shown, e.g.,  in \cite{VP-Tp}. 

Next we wish to show that the external iteration converges. To this aim we work with the solutions to \eqref{eq:double} in iteration $m$ and $m-1$, which we subtract. 
\bas
\Um+ \fF(\Um)+
\tau {A}(\Umm) \Um = \Umm+ \fF(\Umm)+
\tau {A}(\Ummm) \Umm. 
\eas
Introducing $-\tau A(\Umm)\Umm$ on both sides, moving the first two right hand side terms to the left, rearranging and taking inner product with $\Um-\Umm$ as before,  and estimating the monotone terms involving $I+\fF+\tau A(\Umm)$ on the left hand side from below we obtain
\begin{multline}
(1+ \tau \kappa_0) \norm{\Um -\Umm}{}^2  
\\
\leq \tau \inner{{A}(\Ummm) \Umm -\A(\Umm)\Umm ,\Um-\Umm}
\\
\leq 
\tau L_A \norm{\Umm}{} \norm{\Umm - \Ummm}{} \norm{\Um - \Umm}{}.
\end{multline}
Now we can apply the a-priori bound \eqref{eq:uboundA} from \rthe{th:u} which holds for $\Umm$ and we see that the double iteration converges provided 
\ba
\frac{\tau L_A}{(1+\tau \kappa_0)^2} \norm{g}{}<1.
\ea
Even though this is a restrictive bound on $\tau$, at least it does not involve $L_F$. We summarize these calculations below.
\begin{prop}
\label{prop:double}
Under the assumptions of Lemma~\ref{lem:mF} the double iteration \eqref{eq:double} is well-posed and converges. 
\end{prop}

The efficiency of \eqref{eq:double} depends on the ability find a good guess for the Newton step so that the semismooth Newton step converges indeed q-superlinearly.

%%%%%%%%
\subsection{Well-posedness and convergence for non-equilibrium problem \eqref{eq:neqsysn}} 
\label{sec:NEQa}
Now we wish to study \eqref{eq:neqsysn} which we rewrite suppressing $n$, and plugging the expression for $\Chi^n$ into that for $U^n$. In every time step $t_n$ we solve 
\ba
\label{eq:nequn}
U^n  + (1-\bar{B}) \fF(U^n) 
  +\tau A(U^n)U^n= g^n=\tau f^n+U^{n-1}+(1-\bar{B})\Chi^{n-1},
\ea
Suppressing $n$ we see we solve
\ba
\label{eq:nequ}
U  + (1-\bar{B}) \fF(U) 
  +\tau A(U)U=g.
\ea
Since $1>\bar{B}>0$, we see immediately that one can apply the analysis of Section~\ref{sec:EQa}. In particular. we apply the double iteration \eqref{eq:double} and Proposition~\ref{prop:double}.  

%%%%%%%%%%%%%%
\subsection{Well-posedness and convergence for hysteresis problem \eqref{eq:hystsysn}} 
\label{sec:HYSTa}

Now we focus on solving \eqref{eq:hystsysnu} to which we plug \eqref{eq:hystsysnv}. At time step $t_n$, we solve
\ba
\label{eq:hystun}
U^n+
\fF(U^n) + (\resab{0}{\beta_j^n}{\Chi_j^{n-1}-F(U_j^n)})_j
+\tau A(U^n)U^n = g^n
\\
\nonumber 
g^n=\tau f^n + U^{n-1}+\Chi^{n-1}.
\ea
With the sequential approach, $(\beta_j)_j$ is known. In what follows we suppress  
the dependence of $R(0,\beta;v)$ on its first two arguments and write $R(v)$ instead. 

Thus, given a fixed $\Chi,\beta,g$ we must solve for $u$
\bas
U+
\fF(U) + (R(\Chi_j-F(U_j)))_j
+\tau \bar{A}U = g.
\eas
The term $\fF_R(U)=\fF(U) + (R(\Chi_j-F(U_j)))_j$ requires additional analysis. 

We prove that $\fF_R$ is Lipschitz componentwise with constant $L_F$ and monotone. Then  we apply the analysis of Section~\ref{sec:EQa}. 

\begin{lemma}
\label{lem:rf}
Let $F$ be monotone nondecreasing with Lipschitz constant $L_F$, and let $R=R(0,\beta;v)$ be as defined in \eqref{eq:defr}. Then given a fixed $x\in \R$, the function
\ba
\label{eq:deffr}
F_R(x; v) = F(v)+R(x-F(v)).
\ea
is differentiable a.e. in $v$, Lipschitz with constant $L_F$ and monotone in $v$: 
\ba
\label{eq:frmon}
\inner{v_1-v_2,F_R(x,v_1)-F_R(x,v_2)}\geq 0. 
\ea
\end{lemma}
\begin{proof}
The proof is elementary but tedious since $R$ is only differentiable a.e., and $R(x-F(v))$ is anti-monotone in $V$.   

Lipschitz continuity of $F_R$ is immediate; the calculation of $L_{FR}$ follows after we rewrite setting $w_k=F(v_k),y_k=x-w_k$
\begin{multline}
F_R(x;v_1)-F_R(x;v_2)= w_1+R(x-w_1)-(w_2+R(x-w_2))
\\
=y_2-y_1+R(y_1)-R(y_2)=(I-R)(y_2)-(I-R)(y_1).
\end{multline}
Since $I-R$ has Lipschitz constant $1$, thus 
\ba
\norm{F_R(x;v_1)-F_R(x;v_2)}{} = \norm{(I-R)(y_1-y_2)}{} 
\leq \norm{w_1-w_2}{} \leq L_F \norm{v_1-v_2}{}.
\ea

To prove monotonicity, recall $0\leq R'\leq 1$ a.e. Substitute $w_k=F(v_k),\xi_k=x-w_k,k=1,2$, and rewrite
$
F_R(x;v_1)-F_R(x;v_2)=w_1-w_2 + R(x-w_1)-R(x-w_2).
$
For the left hand side of \eqref{eq:frmon} we have
\begin{multline}
\label{eq:FR}
(v_1-v_2) (F_R(x;v_1)-F_R(x;v_2)) \\
= (v_1-v_2)(w_1-w_2) - (v_1-v_2)(R(\xi_2)-R(\xi_1)).
\end{multline}
From monotonicity of $F,R$ and with $\xi_k=x-w_k$ we have
\ba
\label{eq:mon}
v_1-v_2\geq 0 \mathrm{\ iff\  }w_1-w_2 \geq 0\mathrm{\ iff\ } R(\xi_2)-R(\xi_1)\geq 0.
\ea

We study the components of $\inner{v_1-v_2,R(\xi_2)-R(\xi_1)}$ and see that by \eqref{eq:mon} their product is nonnegative, i.e. equal $\abs{(v_1-v_2)(R(\xi_2)-R(\xi_1))}$. We derive thus by repeated application of \eqref{eq:mon}, \eqref{eq:RL}, and the fact that $\xi_2-\xi_1 = w_1-w_2$
%%%
\begin{multline}
(v_1-v_2)(R(\xi_2)-R(\xi_1)) = 
\abs{(v_1-v_2)(R(\xi_2)-R(\xi_1))}
\\
=
\abs{v_1-v_2}\abs{R(\xi_2)-R(\xi_1)} \leq \abs{v_1-v_2}\abs{\xi_2-\xi_1}
\\
=\abs{v_1-v_2}\abs{w_2-w_1}=(v_1-v_2)(w_1-w_2).
\end{multline}
%%%
Applying this identity, \eqref{eq:mon} and \eqref{eq:FR}, upon a sign change we obtain \eqref{eq:frmon} from
%%%
\begin{multline}
(v_1-v_2) (F_R(x;v_1)-F_R(x;v_2)) 
\\
\geq (v_1-v_2)(w_1-w_2) - (v_1-v_2)(w_1-w_2)=0.
\end{multline}
\end{proof} 

With these facts in place, we see that we can easily adapt the proof for the equilibrium case regarding well-posedness of \eqref{eq:hystsysn} and its solvability by the double iteration~\eqref{eq:double}. However, while our problem is still semismooth, compared to (EQ) and (NEQ) we have now an additional source of lack of smoothness due to only piecewise smoothness of $R(0,\beta;\cdot)$. This fact makes Newton solver work harder.

%%%%%%%%%%%%
\subsection{Summary of theoretical results}

We have seen that the Proposition~\ref{prop:double} gives convergence of ``double'' iteration \eqref{eq:double} and handling the double nonlinearity with $\fF(U)$ and $A(U)U$ for the equilibrium problem \eqref{eq:equn}. With small modification, it also applies to the nonequilibrium (NEQ) and hysteresis (HYST) problems \eqref{eq:nequn} and \eqref{eq:hystun}. 

With our unified treatment of the algorithms, the main theoretical challenge is not in the individual features of (NEQ) or (HYST) models but rather still in the double nonlinearity. 

We implement the double iteration \eqref{eq:double} in practice for all the variants (EQ), (NEQ), (HYST). The results and examples are presented below.  

%%%%%%%%%%%%
\section{Examples} 
\label{sec:examples}

In this section we show numerical examples for \eqref{eq:uchi} with (EQ), (NEQ), and (HYST) closure. We start with a detailed illustration of calibration of generalized play  (HYST) model in Section~\ref{sec:hysteresis}. 

For this model we illustrate how to work with an ODE system, \mpa{a reduction of \eqref{eq:hystsysn} from $(U^n, \Chi^n) \in R^M \times \R^M$ to $(U^n, \Chi^n) \in \R \times \R$}. We are able to confirm first order convergence  of the numerical scheme; see Section~\ref{sec:ODE}. 

Section~\ref{sec:results-pde} is devoted to the PDE examples and simulations of \eqref{eq:pdeu} coupled to $\chi$ given by the  (EQ), (NEQ), (HYST) variants.

%%%
\subsection{Calibration of generalized play (HYST) model}
Here we work with the notation of the physical model denoting temperature by $\theta$. We aim to develop a hysteresis model for $\chi(t)$ with which $\fr\leq \chi(t)\leq \fl(t)$, where the lower bound is the same as the curve $\fr(t)$ ``on the right'' with subscript $r$. Analogously, $\fl$ represents the bound $\fl(t)$ on the left.  Both $\fl,\fr$ should be at least Lipschitz continuous.

First we discuss in practice how to get $\fl$ and $\fr$. Typically, assume we are given some equilibrium curve $\chi=F(\theta)$. We  aim to define some $G(\theta)$ so that
\ba
\label{eq:FG}
F(\theta) \leq \chi \leq G(\theta),
\ea
where $G$ is similar to $F$ but has an off-set which models the different rate and delay of nucleation. 
By design, we want $G$ and $F$ to agree for $\theta>0$ as well as below certain low temperature $\theta<\theta_0$
\bsa
\label{eq:FGrules}
\ba
G(\theta)=F(\theta), \theta<\theta_0; \; G(\theta)=F(\theta), \theta>0.
\ea
We also want these to be smooth, and we parametrize $G(\theta)$ requiring
\ba
\label{eq:FG0}
F(\theta_0)=G(\theta_0); F'(\theta_0)=G'(\theta_0);
F(0)=G(0).
\ea
\esa

\begin{figure}
\centering
\includegraphics[width=0.45\textwidth]{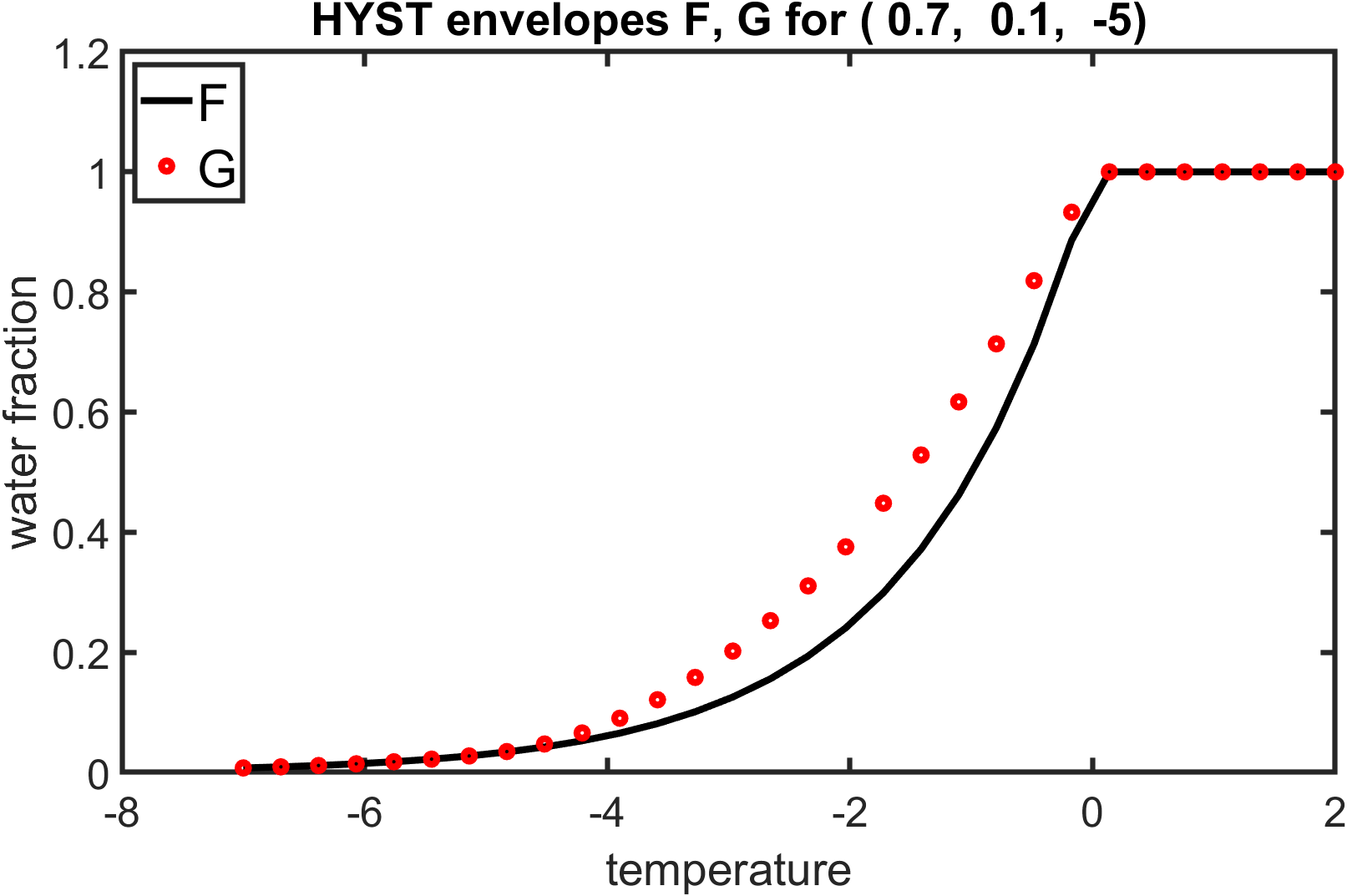}
\includegraphics[width=0.45\textwidth]{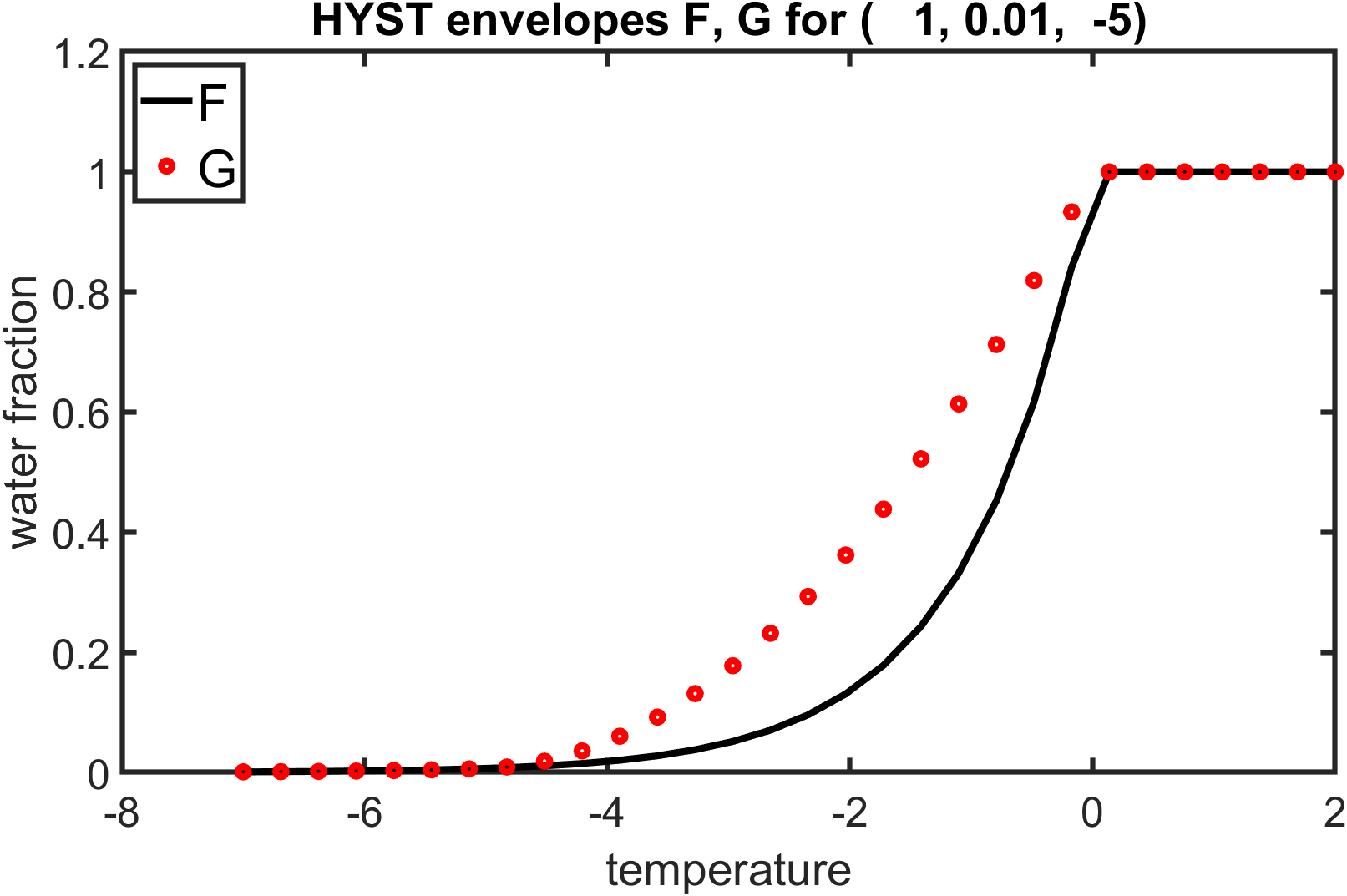}
\caption{The envelopes of the generalized play HYST model calibrated in Example~\ref{ex:hyst-FG} (i) and (ii). The titles indicate the parameters ($b,\bbar,\theta_0$). }
\label{fig:chi-hyst}
\end{figure}

\begin{example}[Calibration of $F(\theta), G(\theta)$ for permafrost application]
\label{ex:hyst-FG}
Let $b>0$ be given and $F(\theta) = e^{b \theta}$ be as in \eqref{eq:chiM}. For some given $\bar{b}>0$ amd $\theta_0<0$ we  propose 
\ba
\label{eq:gdef}
G(\theta) = a e^{\bar{b} \theta} + D\theta + C,
\ea
where $a,D,C$ are found with the three conditions in \eqref{eq:FG0}. We find
\ba \label{eq:hys_values}
    a = \frac{e^{b \theta_0} - b \theta_0 e^{b \theta_0} - 1}{e^{\bar{b} \theta_0} - \bar{b} \theta_0 e^{\bar{b} \theta_0} - 1};
    C = 1 - a ;
    D = b e^{b \theta_0}- a \bar{b} e^{\bar{b} \theta_0}.
\ea
%%%
In particular, we obtain data in Table~\ref{tab:hyst-FG}
%%%
\begin{table}[ht]
\label{tab:hyst-FG}
\caption{Calibration data in Example~\ref{ex:hyst-FG}. }
\centering
\begin{tabular}{l||l|l|l|lll}
\hline
CASE&$b$&$\bbar$&$\theta_0$&$a$&$D$&$C$\\
\hline
(i)&0.7&0.1&-5
%a = 9.579527234657450;
%D = -0.559889528941569;
%C = -8.579527234657450.
&
9.5795&
-0.5598&
-8.5795.
\\
(ii)&1&0.01&-5&
%A =7.936224678392653e+02
%D =-7.542432486536721
%C =-7.926224678392653e+02
%%
793.62
&
-7.5424
&
 -792.6225
\\
\mpc{(iii)}&\mpc{0.5}& \mpc{0.75}&\mpc{-5}&
\mpc{2.3269}&&\mpc{0.0274}
\\
\hline
\end{tabular}
\end{table}
with the graphs $F,G$ plotted in Figure~\ref{fig:chi-hyst}.  
\\
We also consider a different case (iii) when the graph $F$ is shifted to the left, and we impose matching with $G(\theta) = a e^{\bar{b}\theta} + C$ where we find 
$
  a = \frac{b e^{b \theta_0}}{\bar{b} e^{\bar{b} \theta_0}}; 
    C = e^{b \theta_0} - \frac{b}{\bar{b}} e^{b \theta_0}
$. 
\end{example}

%%%

%%%
\subsection{Simulation and convergence of (HYST) for an ODE model}
\label{sec:ODE}
The envelopes for HYST model developed in Example~\ref{ex:hyst-FG} can now be applied in the evolution model \eqref{eq:hystFG}. We set-up an example to illustrate that the generalized play model works well and delivers solutions which stay within these envelopes and satisfy \eqref{eq:FG}. 

\begin{example}
\label{ex:hyst-ODE}
We calibrate the hysteresis data $(0,\beta)$ 
given $F(u)$ and $G(u)$ as in Example~\ref{ex:hyst-FG}(ii) and (iii). We solve the ODE \eqref{eq:hystFG} with scheme \eqref{eq:hystFGn} applying to the resolvent of $\cg(0,\beta;\cdot)$ the sequential strategy for $\beta(t)=G(u(t))-F(u(t))$ based on an assumed input $u(t)$ 
\ba
u(t) = h(t) \cos \left (\frac{\pi}{4} t \right ) + g(t), \;
    h(t) = \left \{ \begin{array}{cc}
        8 & t < 4 \\
        4 & t \geq 4
    \end{array}\right., \;
    g(t) = \left \{ \begin{array}{cc}
        -2 & t < 4 \\
        \frac{t}{2} - 8 & t \geq 4.
    \end{array}\right. 
\ea
We consider the interval $(0, T], T=30$ with a timestep of $\tau=3.75 \times 10^{-2}$.  
\end{example}

\begin{figure}
\centering
\includegraphics[width=0.45\textwidth]{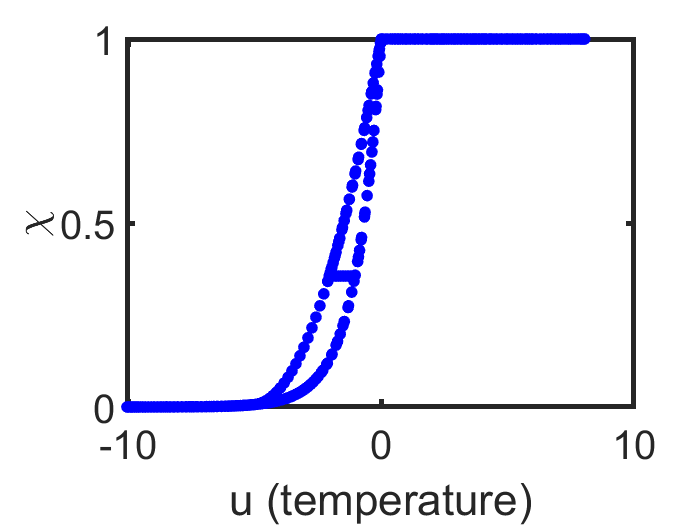}
\includegraphics[width=0.45\textwidth]{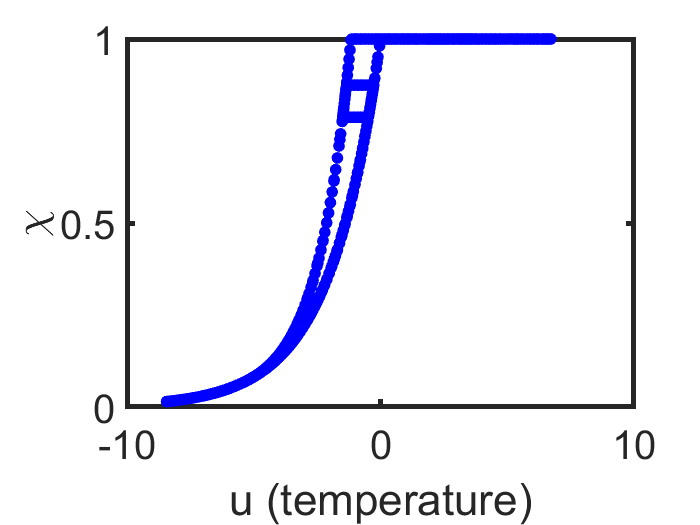}
\caption{Phase plot of solutions $(U^n,\Chi^n)_n$ to Example~\ref{ex:hyst-ODE} case (ii) and (iii).}
\label{fig:hyst-ODE}
\end{figure}

Figure~\ref{fig:hyst-ODE} shows computed solution which adheres well to the prescribed envelopes and travels horizontally between these envelopes $F(u)$ and $G(u)$. 
%%%

%%%%%%%%%%%
\subsubsection{Coupled ODE with generalized play (HYST) model}
\label{sec:HYST-ode}

Next we consider a coupled model where $u(t)$ is not explicitly given as in Example~\ref{ex:hyst-ODE} but rather is governed by its own dynamics
\ba \label{eq:theta_dynamics}
    \ddt u + \ddt \chi + A u = f(t); \; t \in (0,T]; u(0)=\uz, \chi(0)=\chi_{init},
\ea
coupled to \eqref{eq:hystFG}. 
Here $A > 0$ is fixed and $f$ is some input function. We use Newton solver \eqref{eq:newton} to solve the discrete step \eqref{eq:hystFGn}. 

\begin{example}
\label{ex:hyst-ODEsys}
We let $A = 0.02$ and design an oscillating forcing function 
\ba
    f(t) = h(t) \cos \left (\pi t \right ) + g(t); \;
    h(t) = \left \{ \begin{array}{cc}
        16 & t < 1 \\
        4 & t \geq 1
    \end{array}\right. ;\;
    g(t) = \left \{ \begin{array}{cc}
        -15 & t < 1 \\
        4t - 30 & t \geq 1.
    \end{array}\right. 
\ea
The hysteresis graph is calibrated with $b = 1$, $\bbar = 0.1$, $\theta_0 = -5$. We use $\uz = -0.2$, $\chi_{init} = e^{-0.5}$, and consider the time interval $(0,T], T= 10$. 
\end{example}

For this example we study the convergence of the solutions as $\tau \to 0$. Since we do not know the true solution, we use as a proxy the fine grid solution with $\tau = 0.0001$.  
We compute the order of convergence for the solution vector $\zeta := (\chi, u)$ considering for each $\tau$, $\zeta_{\tau}(t)-\zeta_{fine}(t)$. In  Table~\ref{tab:order_of_convergence} we report the errors $\norm{\zeta_{\tau}-\zeta_{fine}}{}$, as well as the order of convergence. We see first order convergence in all norms.

\begin{table}[ht]
  \caption{Error and convergence order for Example~\ref{ex:hyst-ODE}.}
    \label{tab:order_of_convergence}
    \centering
    \begin{tabular}{ccccccc}
        \hline
        $\tau$ & $||\zeta_{\text{fine}}- \zeta_\tau||_1$ & Order & $||\zeta_{\text{fine}} - \zeta_\tau||_2$ & Order & $||\zeta_{\text{fine}} - \zeta_\tau||_\infty$ & Order \\
        \hline
        %%%%
        0.1 & 0.1750 & - & 0.01750 & - & 0.1750 & - \\
        0.01 & 0.0149 & 1.0702 & 0.0149 & 1.0702 & 0.0149 & 1.0702\\
        0.001 & 0.0013 & 1.0495 & 0.0013 & 1.0495 & 0.0013 & 1.0495\\
        %%%%
        \hline
    \end{tabular}
\end{table}

%%%
\subsection{Computational results for PDE model}
\label{sec:results-pde}
Now we show simulation results for the PDE model \eqref{eq:pdeu} with each one of \eqref{eq:EQ}, \eqref{eq:NEQ}, \eqref{eq:HYST} simulated with the schemes \eqref{eq:eqsysn}, \eqref{eq:neqsysn}, \eqref{eq:hystsysn}, respectively. We use Newton's method and the $A$-lagging approach \eqref{eq:double} for resolving the double nonlinearity. 

The examples are designed to show robustness of the scheme, with $(u(x,t),\chi(x,t))$ responding to the varying boundary conditions. Such conditions are realistic in permafrost examples where the surface boundary temperature varies daily.

We first focus on the (EQ) model, and next consider (NEQ) and (HYST) cases. 

{\bf Data for \eqref{eq:pdeu}:} We set $F(u)=e^{b u}$ with $b=1$. We also premultiply \eqref{eq:pdeall} by $10^6$, and set the time scale of $10^6\mpunit{sec} \approx 11.57~\mpunit{day}$. To calculate $c(u),k(u)$ we assume $\eta=0.32$ and use data from \cite{BPV} in Example~\ref{ex:ck}.  Now we use  $\ttc=\tc$, and $\ttk=10^6\tk$, and  $\ttk_f=2.06 \times 10^{-2},\ttk_u=1.51 \times 10^{-2}$, $\ttc_f = 2.21\times 10^{-2}; \ttc_u = 2.94\times 10^{-2}$. We have
%%
%% apprximate values when eta=0.32 (multiply by 1e-2)
%% cf = 2.21; cu = 2.94
% %% kf = 2.06; ku = 1.511
%%
\bas
k(u)&=&\ttk_f+(\ttk_u-\ttk_f)F(u);
\\
c(u)&=& \begin{cases}(\ttc_u-\ttc_f) \int_0^u F(v)dv+\ttc_f u, &(u<0), \\\ttc_u u, &(u\geq 0).\end{cases}
\eas
\begin{example}
\label{ex:pde-EQ}
Let $\Omega=(0,1)$, and $T=3$ (about 33$~\mpunit{days}$).  We set the initial and boundary conditions 
\bsa
\ba
\uz^{EQ}(x)=-5,
%\\
\label{eq:pde-bc}
u(0,t)=\begin{cases}5, &t\leq 1; 
\\
-10(t-1)+5, &1<t\leq 2;
\\
10(t-2)-5, &2 <t, 
\end{cases}
\;\;\; u(1,t)=-5.
\ea
\esa
%fd1d_heat_hysteresis(50,0.1,15,0,0);
%vd = -2*(t-5)+5; vup = 2*(t-10)-5;

%%
We conduct simulations with discretization parameters $M=100$, $\tau=0.01$, considering also testing with other $M,\tau$. 
\end{example}

\begin{figure}[ht]
\centering
\includegraphics[width = 0.3\textwidth]{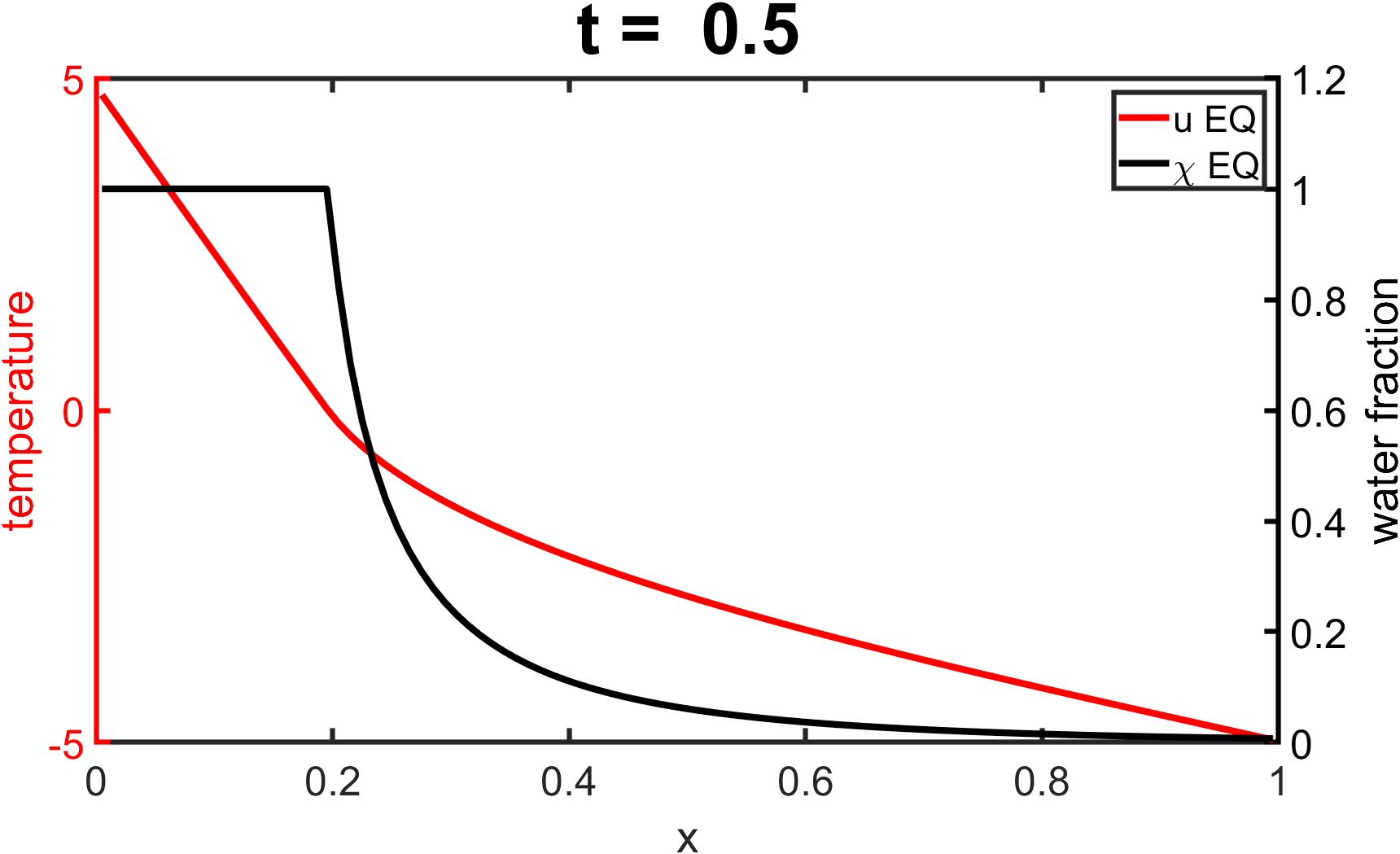}
\includegraphics[width = 0.3\textwidth]{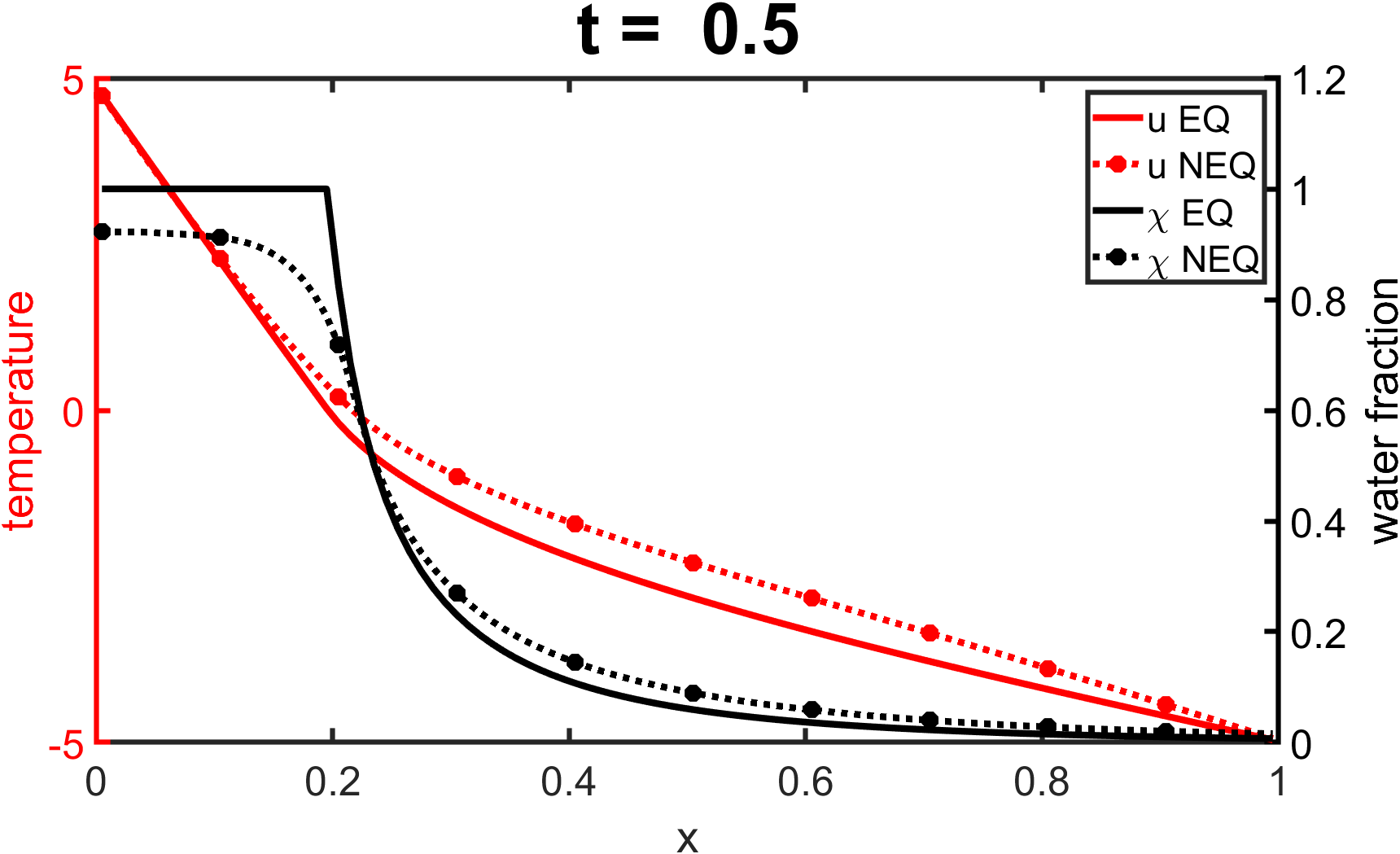}
\includegraphics[width = 0.3\textwidth]{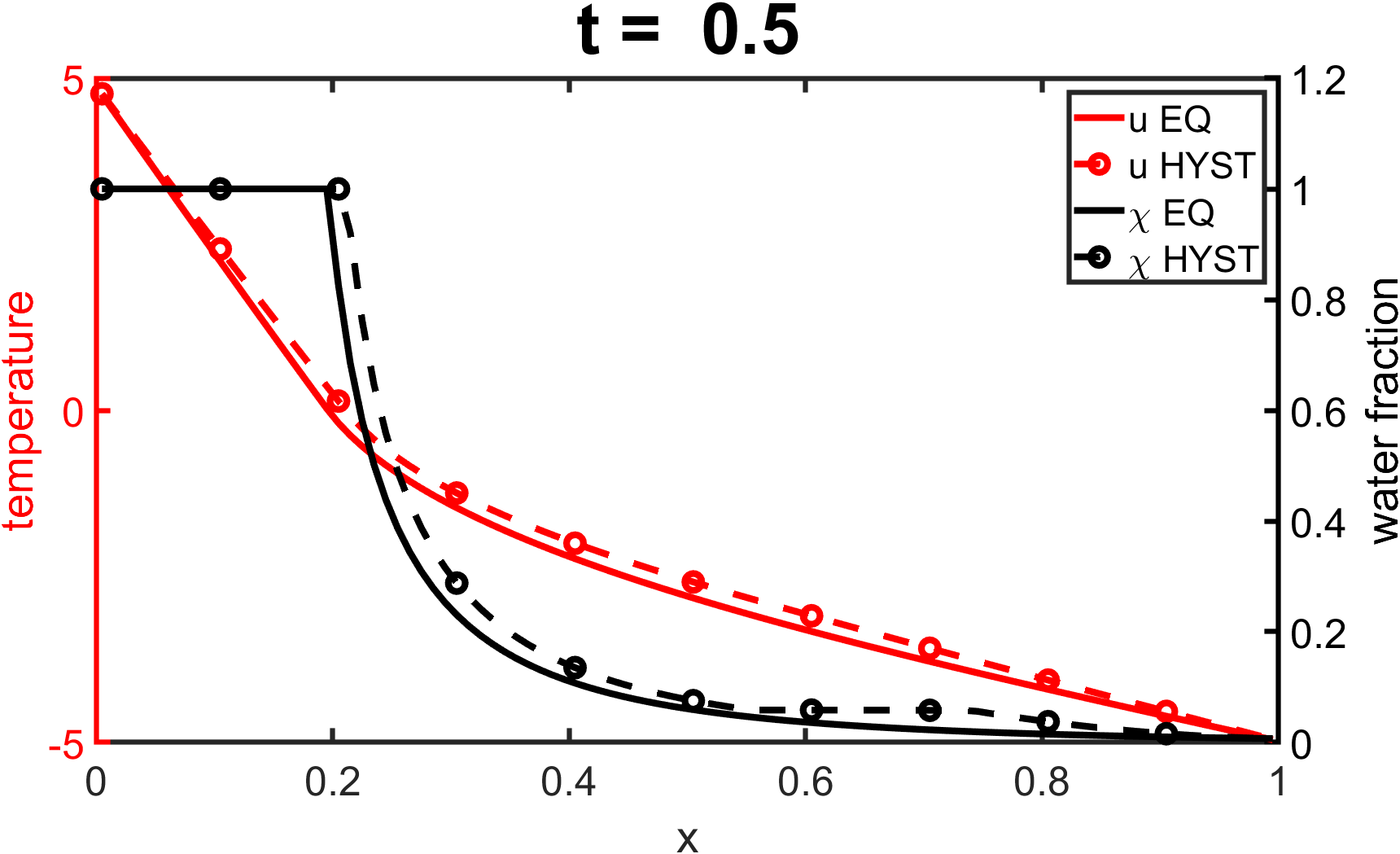}
\\
\includegraphics[width = 0.3\textwidth]{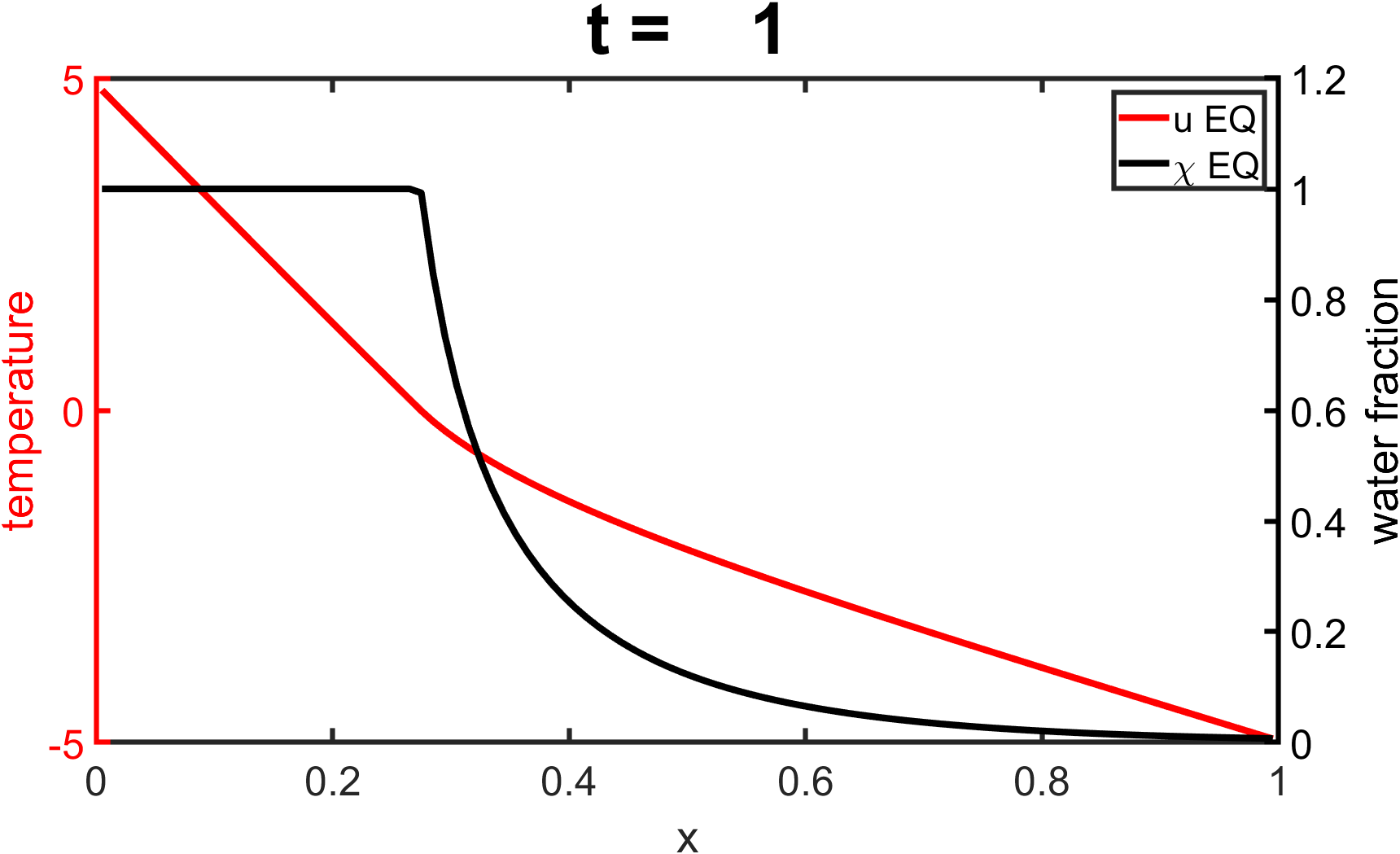}
\includegraphics[width = 0.3\textwidth]{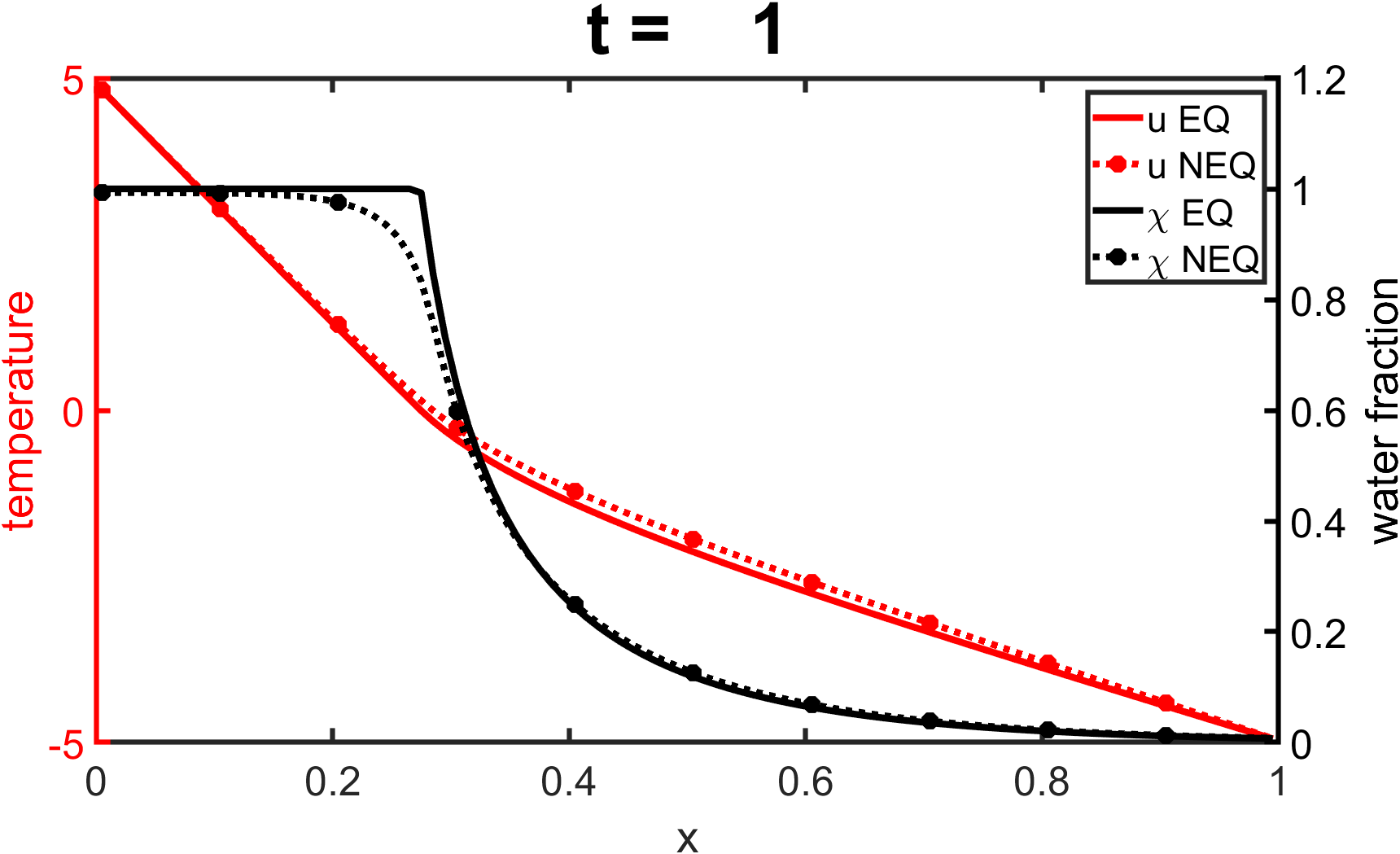}
\includegraphics[width = 0.3\textwidth]{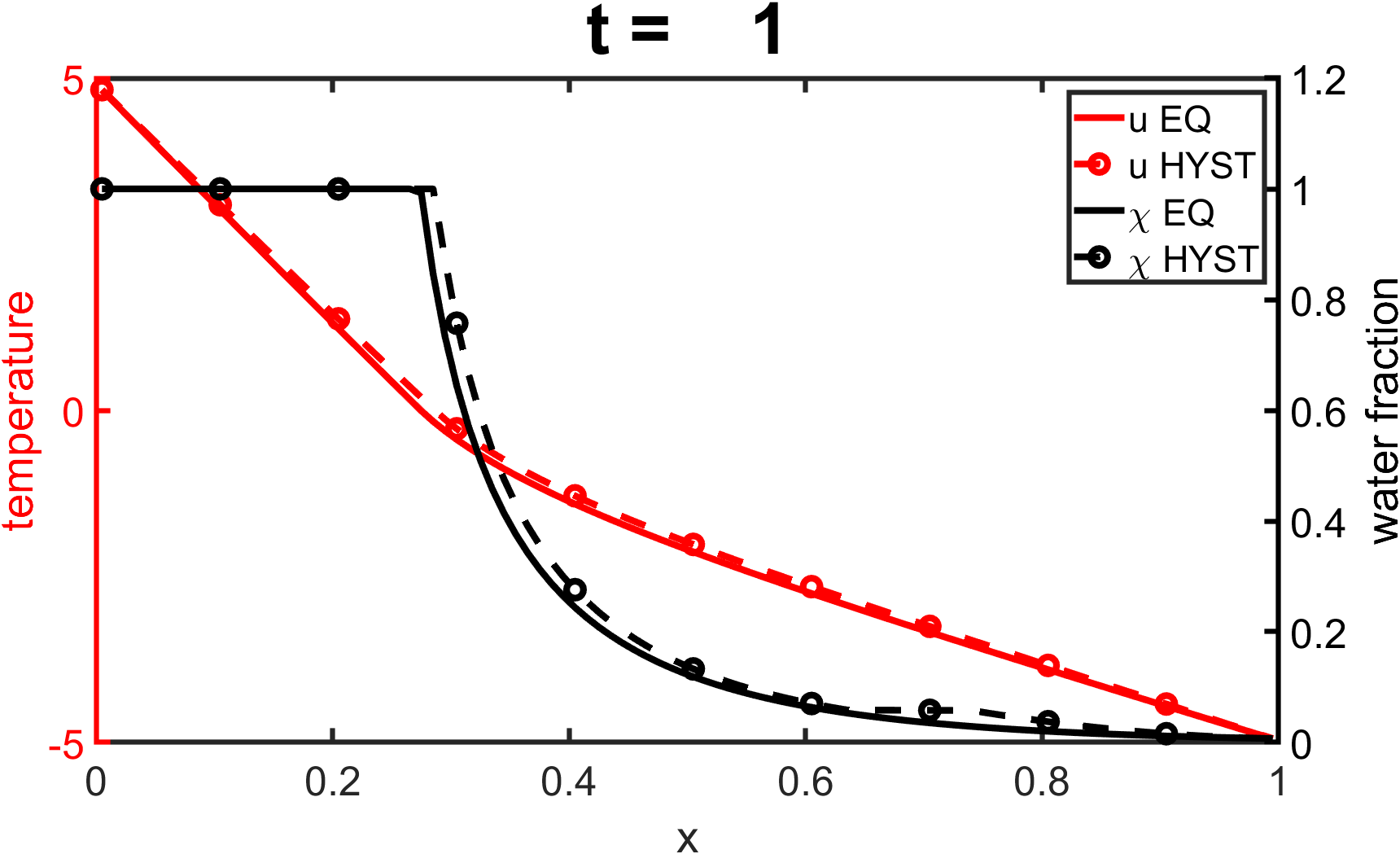}
\\
\includegraphics[width = 0.3\textwidth]{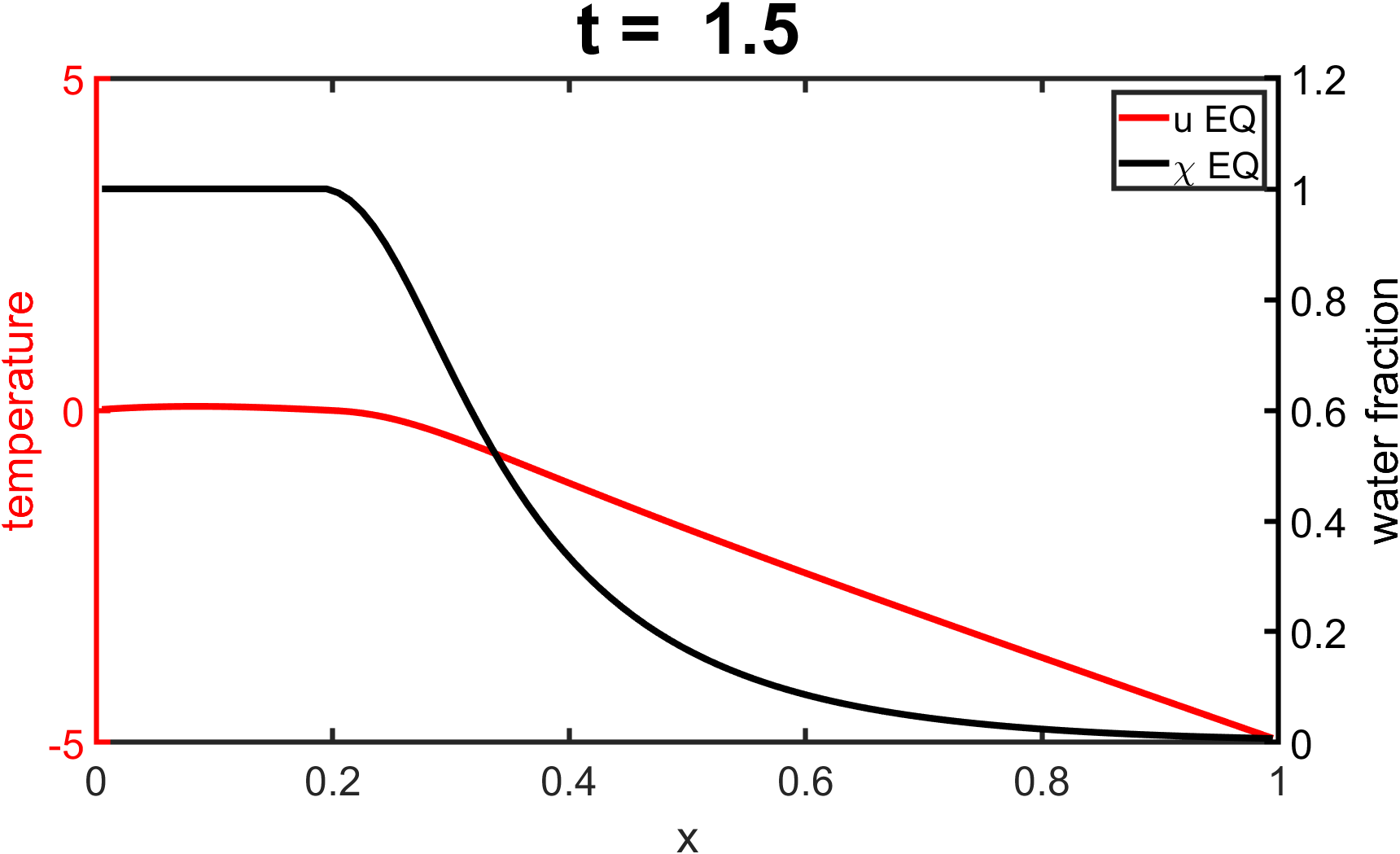}
\includegraphics[width = 0.3\textwidth]{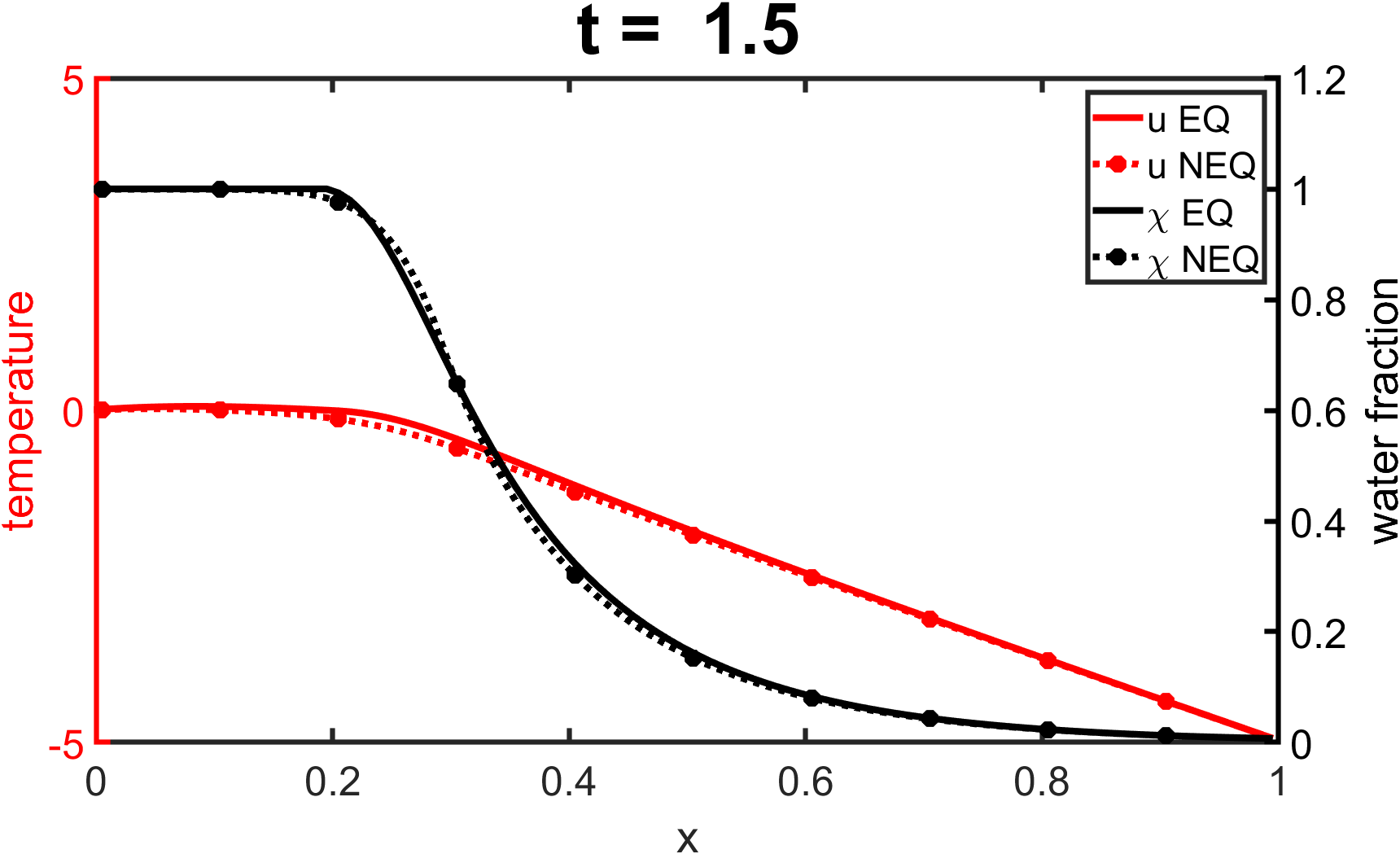}
\includegraphics[width = 0.3\textwidth]{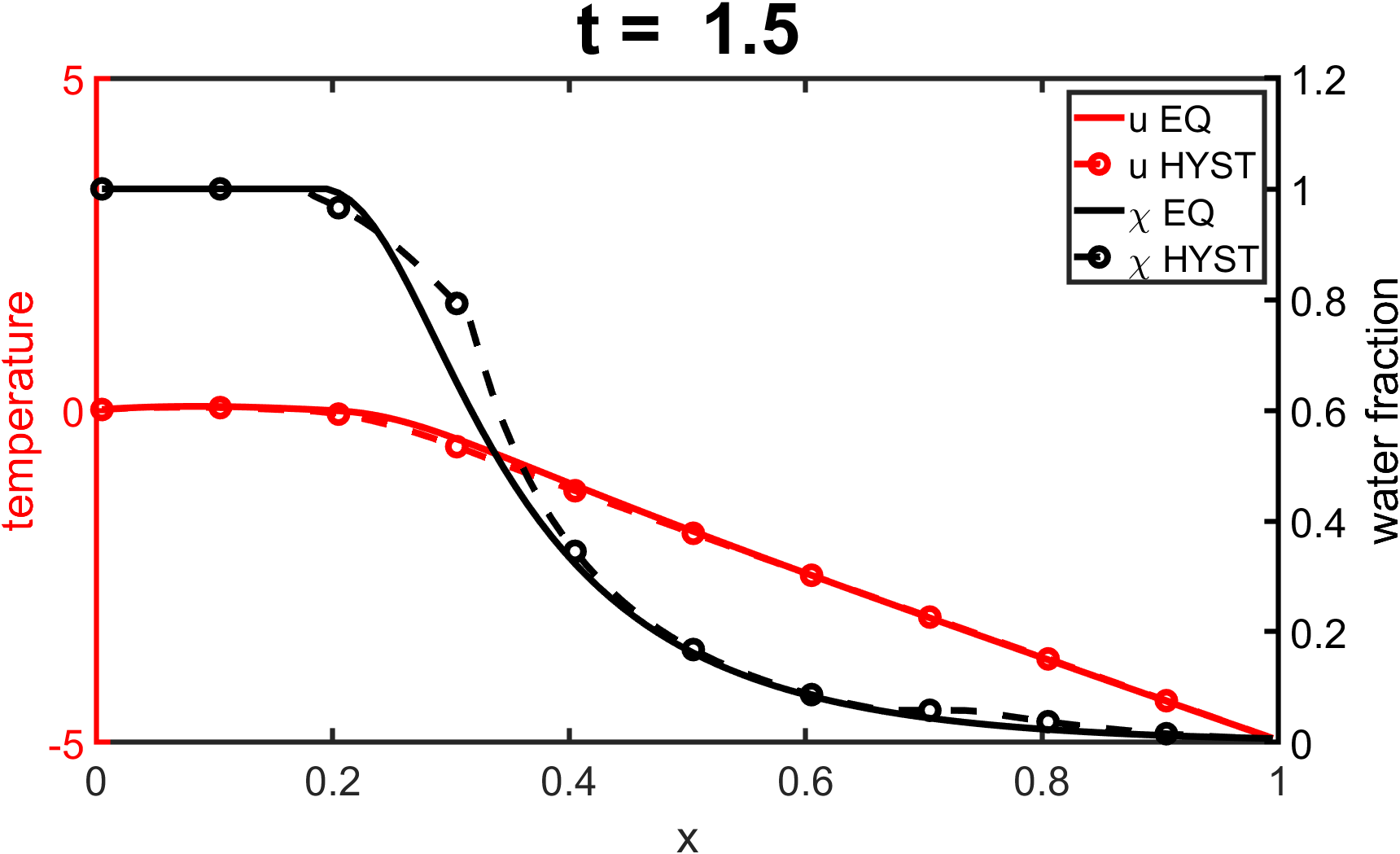}
\\
\includegraphics[width = 0.3\textwidth]{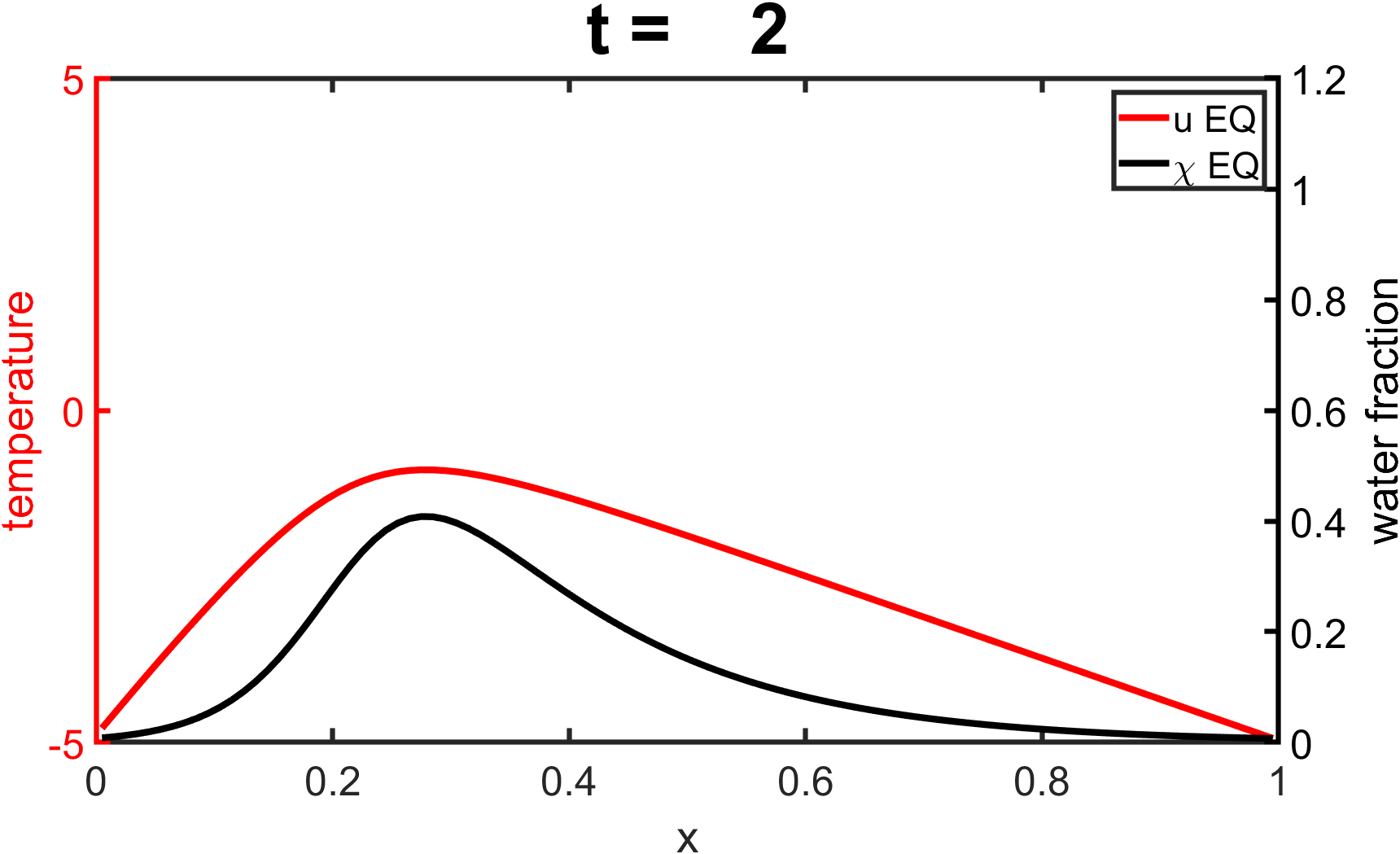}
\includegraphics[width = 0.3\textwidth]{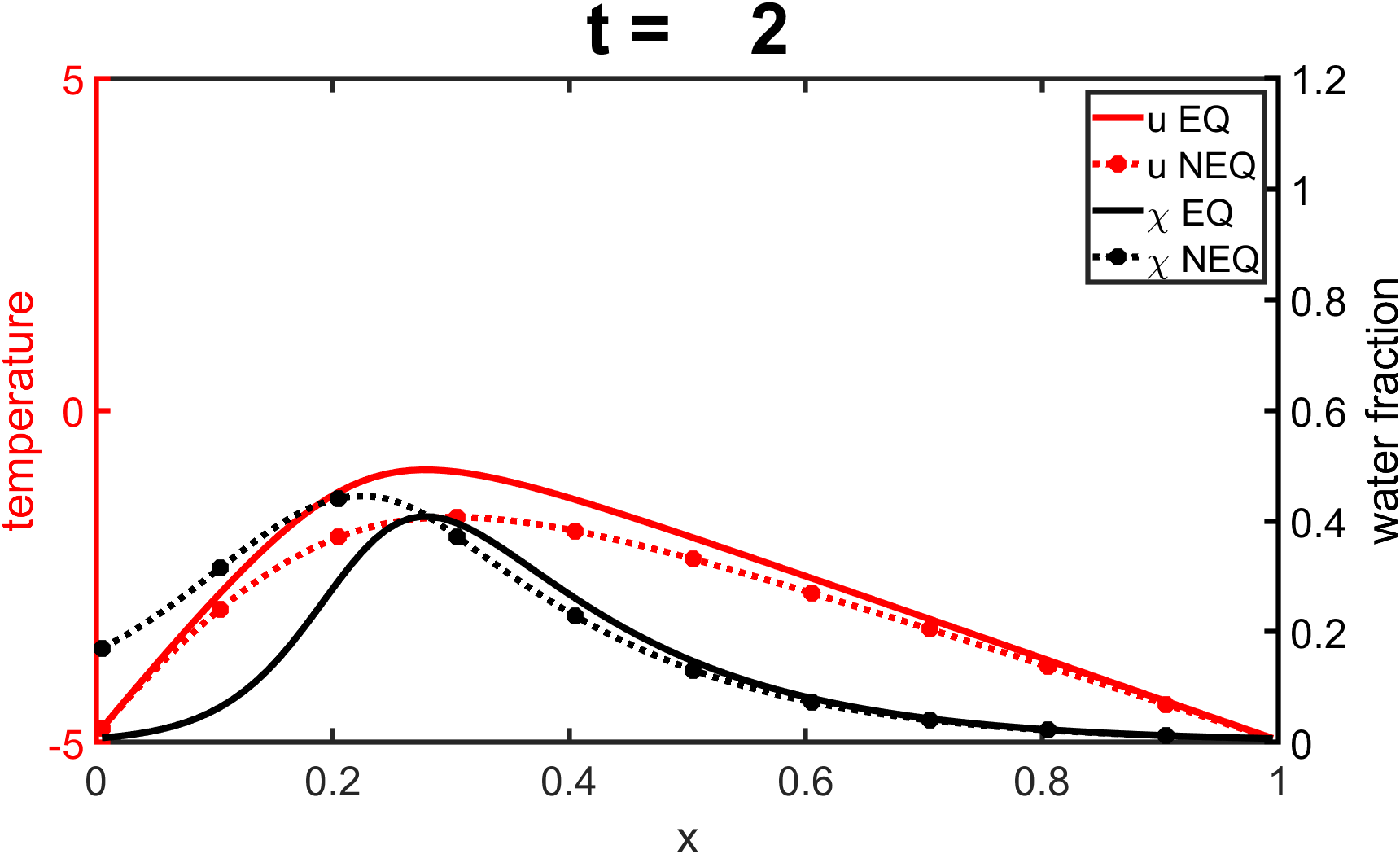}
\includegraphics[width = 0.3\textwidth]{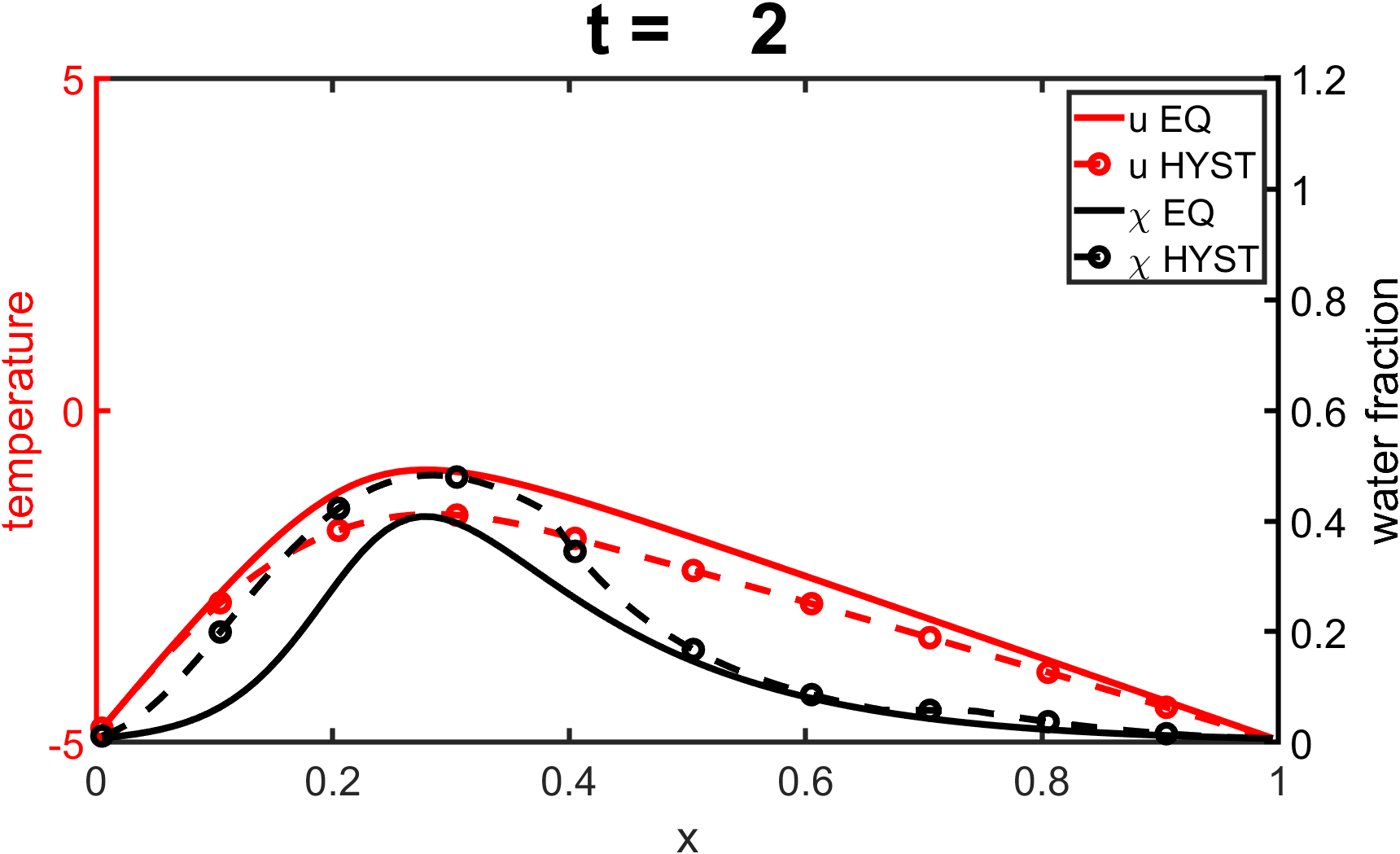}
\\
\includegraphics[width = 0.3\textwidth]{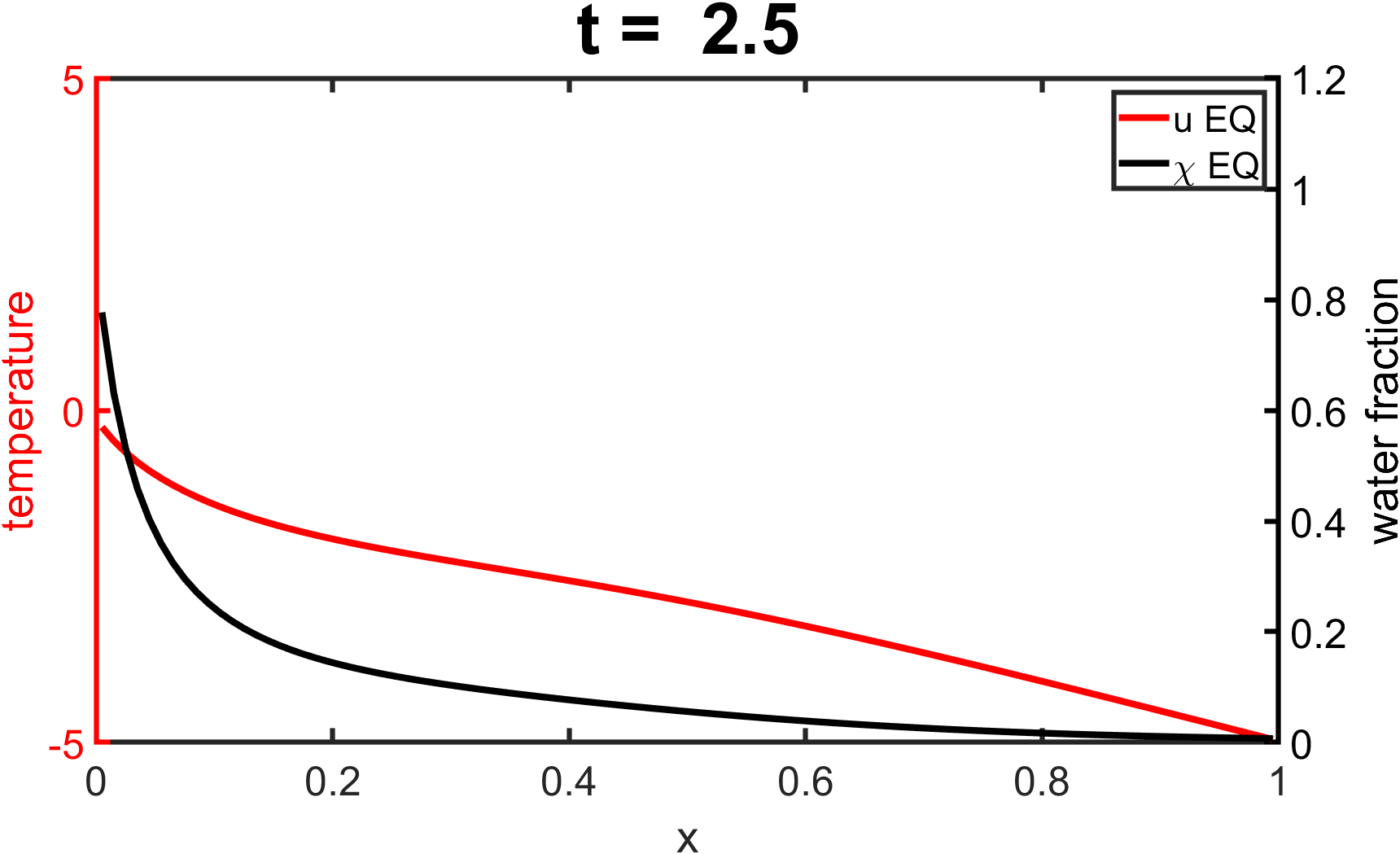}
\includegraphics[width = 0.3\textwidth]{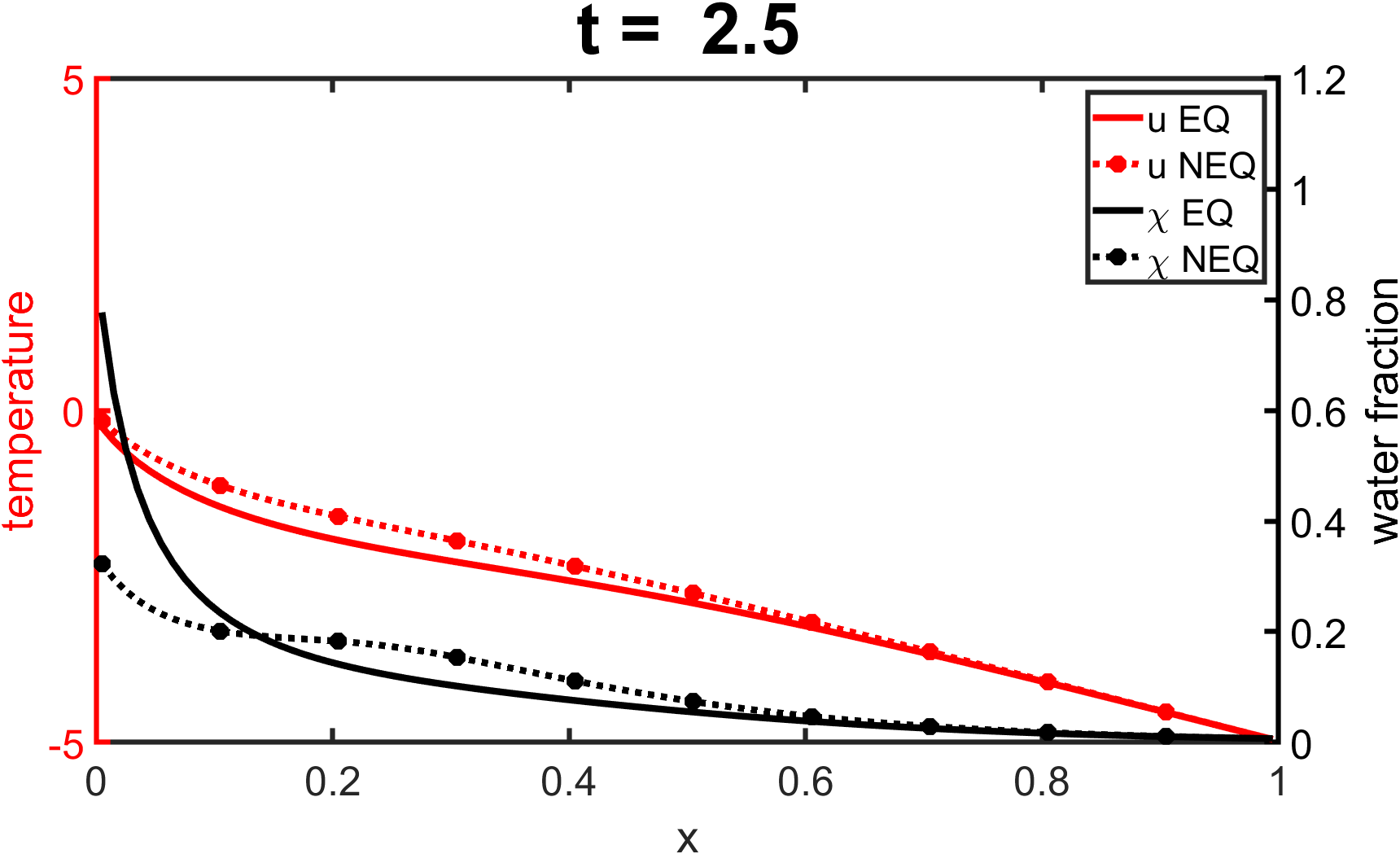}
\includegraphics[width = 0.3\textwidth]{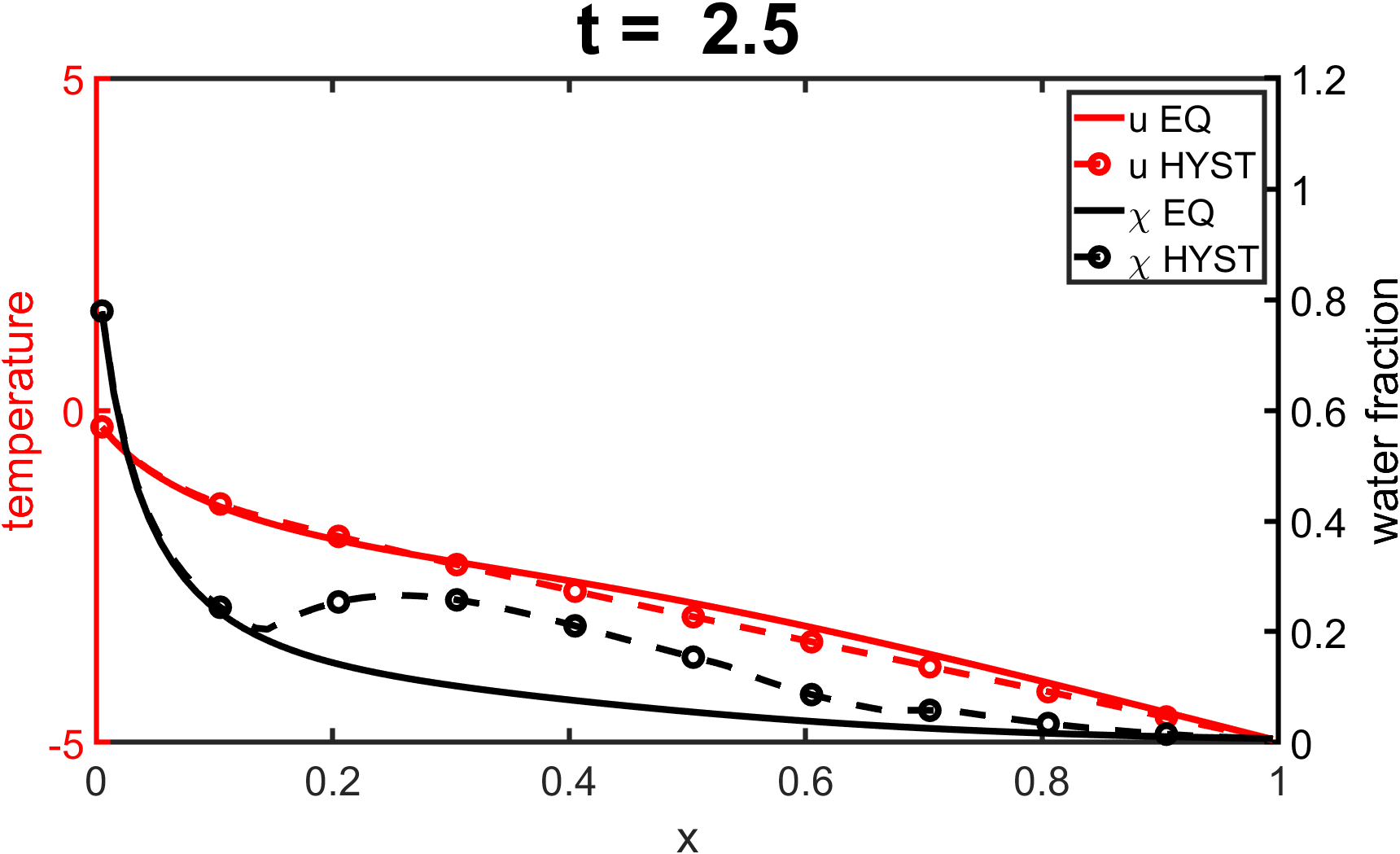}
\\
\includegraphics[width = 0.3\textwidth]{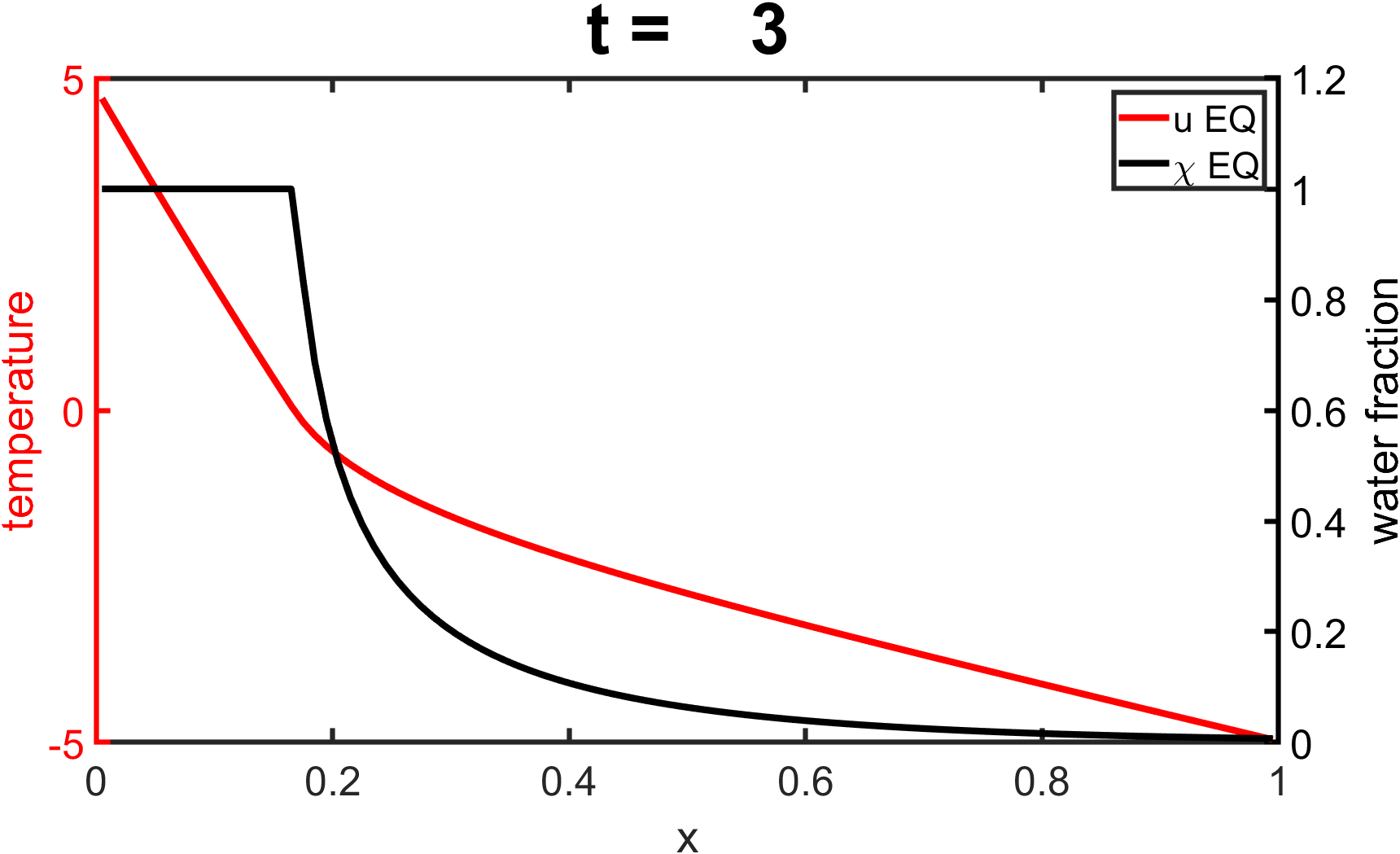}
\includegraphics[width = 0.3\textwidth]{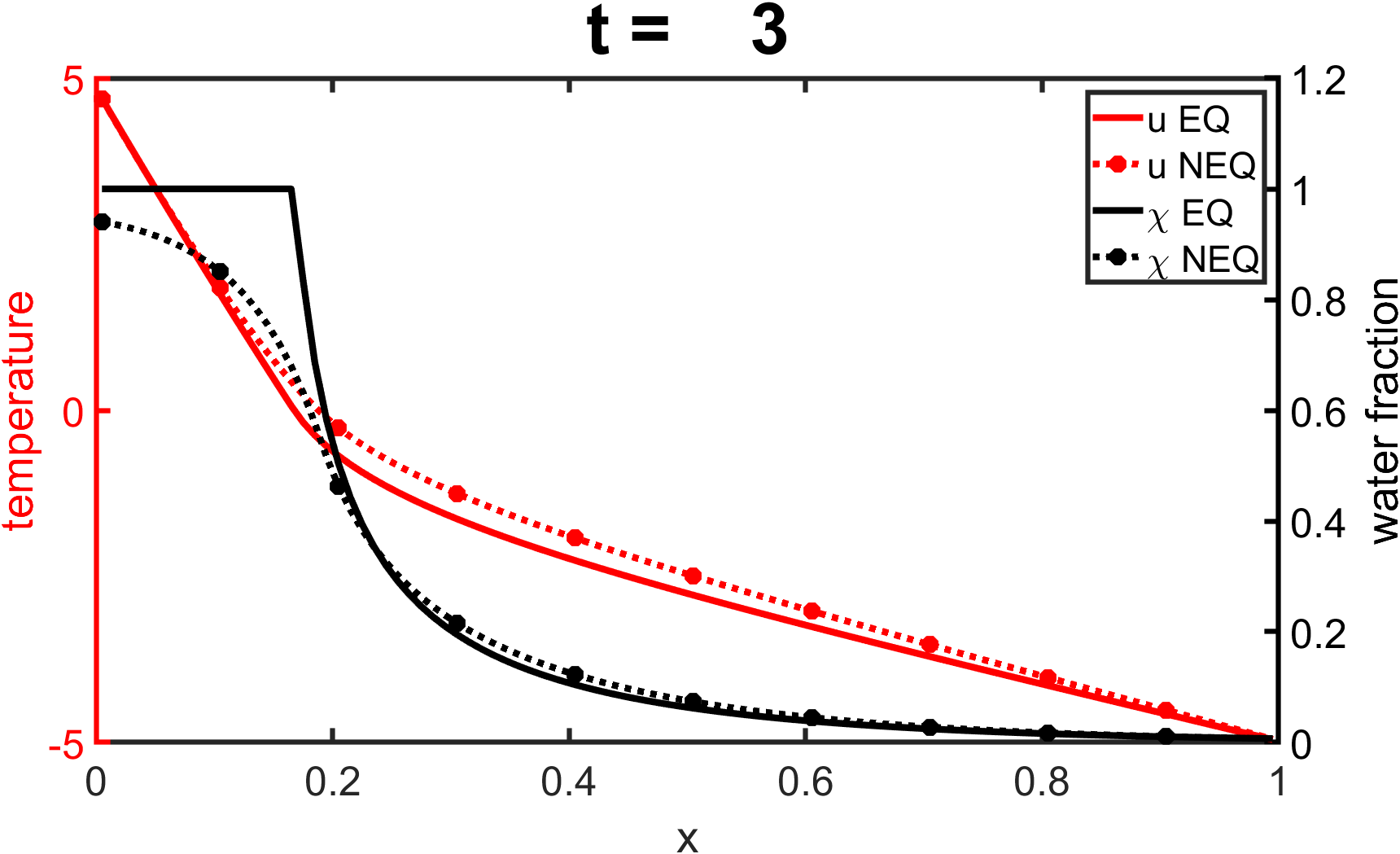}
\includegraphics[width = 0.3\textwidth]{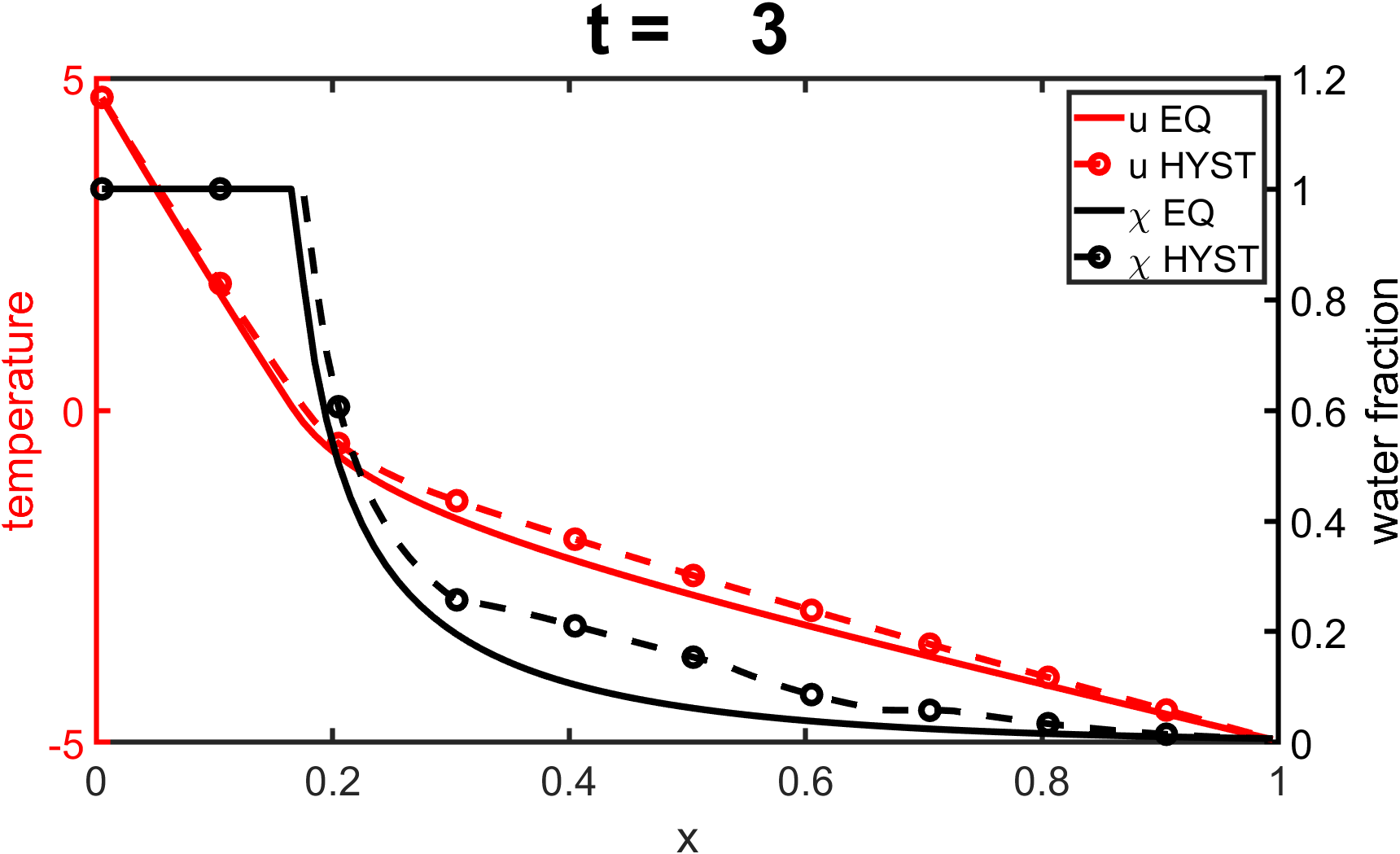}
\caption{Solutions to Example~\ref{ex:pde-EQ} and \ref{ex:pde-ALL} at time steps $t=0.5,1,1.5,2,2.5,3$. We plot both the temperature $u(x,t)$ (scale on the left axis) and the water fraction  $\chi(x,t)$ (scale on the right axis). In the left column we plot only the solutions to the (EQ) model, as indicated in the legend. In the middle and right we plot the solutions to the (NEQ) and (HYST) models, respectively, with the equilibrium solution (EQ) plotted for reference.}
\label{fig:pde}
\end{figure}

We present the plots of evolution of $u(x,t),\chi(x,t)$ in Figure~\ref{fig:pde} (left), with each row corresponding to selected time steps 
$t=0.5,1,1.5,2,2.5,3$.

We see how the temperature $u(x,t)$ adapts to the evolving boundary condition $u(0,t)$, with the water fraction following suit, and the free boundary apparent from the plot of $\chi(x,t)$ at the points $x$ where $u(x,t) \approx 0$. In particular, the temperature $u(x,t),t\leq 5$ increases gradually from $-5$ to $5$ while maintaining the right boundary condition~\eqref{eq:pde-bc}. When $t>1$, $u(x,t)$ responds to the fluctuating upper boundary condition $u\vert_{x=0,t}$.  We also see good resolution of the free boundary; see, e.g., the plot of $\chi(x,t)$ at $t=0.5$ (in the image in the upper left corner). 

\begin{example}
\label{ex:pde-ALL}
We use the same data as in Example~\ref{ex:pde-EQ} for the (NEQ) and (HYST) models. These start out of equilibrium with 
\ba
\chi^{NEQ}_{init}(x)=F(\uz(x))+0.1,
\\
\chi^{HYST}_{init}(x)=F(\uz(x))+0.1.
\ea
In the (NEQ) model we use $B=5$.  We also use $B=10$ and $B=20$ for comparison.
\\
For (HYST) model we calibrate the hysteresis model as shown in Example~\ref{ex:hyst-FG}(ii). 
\end{example}

We simulate $(u^{NEQ}(x,t),\chi^{NEQ}(x,t))$
and $(u^{HYST}(x,t),\chi^{HYST}(x,t))$. 
We add the plots of their approximations to Figure~\ref{fig:pde} to supplement the equilibrium model plotted for reference. 

We see that the solutions $\chi^{NEQ},\chi^{HYST}$ are  distinct from $\chi^{EQ}$ reaching even close to $\norm{\chi^{EQ}-\chi^{NEQ}}{\infty}\approx 0.85$, with the corresponding difference for the temperature $\norm{u^{EQ}-u^{NEQ}}{\infty}\approx 1$. Similarly, we have $\norm{\chi^{EQ}-\chi^{HYST}}{\infty}\approx 0.85$, $\norm{u^{EQ}-u^{NEQ}}{\infty}\approx 1.3$.

The solutions $\chi^{NEQ}(x,t)$ ``follow'' $\chi^{EQ}(x,t)$ with some delay proportional to $1-\bar{B}$. Since the free boundary \mpc{is} present and both $\chi,u$ are evolving, \mpc{the variable $\chi^{NEQ}(x,t)$ is not always very close} to $\chi^{EQ}(x,t)$. The simulations capture the richness of the dynamics due to  variable boundary data. 

In turn, the profile $\chi^{HYST}(x,t)$ adheres to the constraints in \eqref{eq:hystFG} which only become active when the local values of $u(x,t)$ change from increasing to decreasing. In the end $\chi^{HYST}(x,t)$  follows, as expected,  either $F(u(x,t))$ or $G(u(x,t))$, or remains between these, even if it is not easy to track down exactly where the constraint is active.  

The solutions corresponding to the (NEQ) case are smoother than those for (HYST). However, of course both can be expected to be only as smooth as the (EQ) solutions. 

\begin{figure}[ht]
\centering
\includegraphics[width=.45\linewidth]{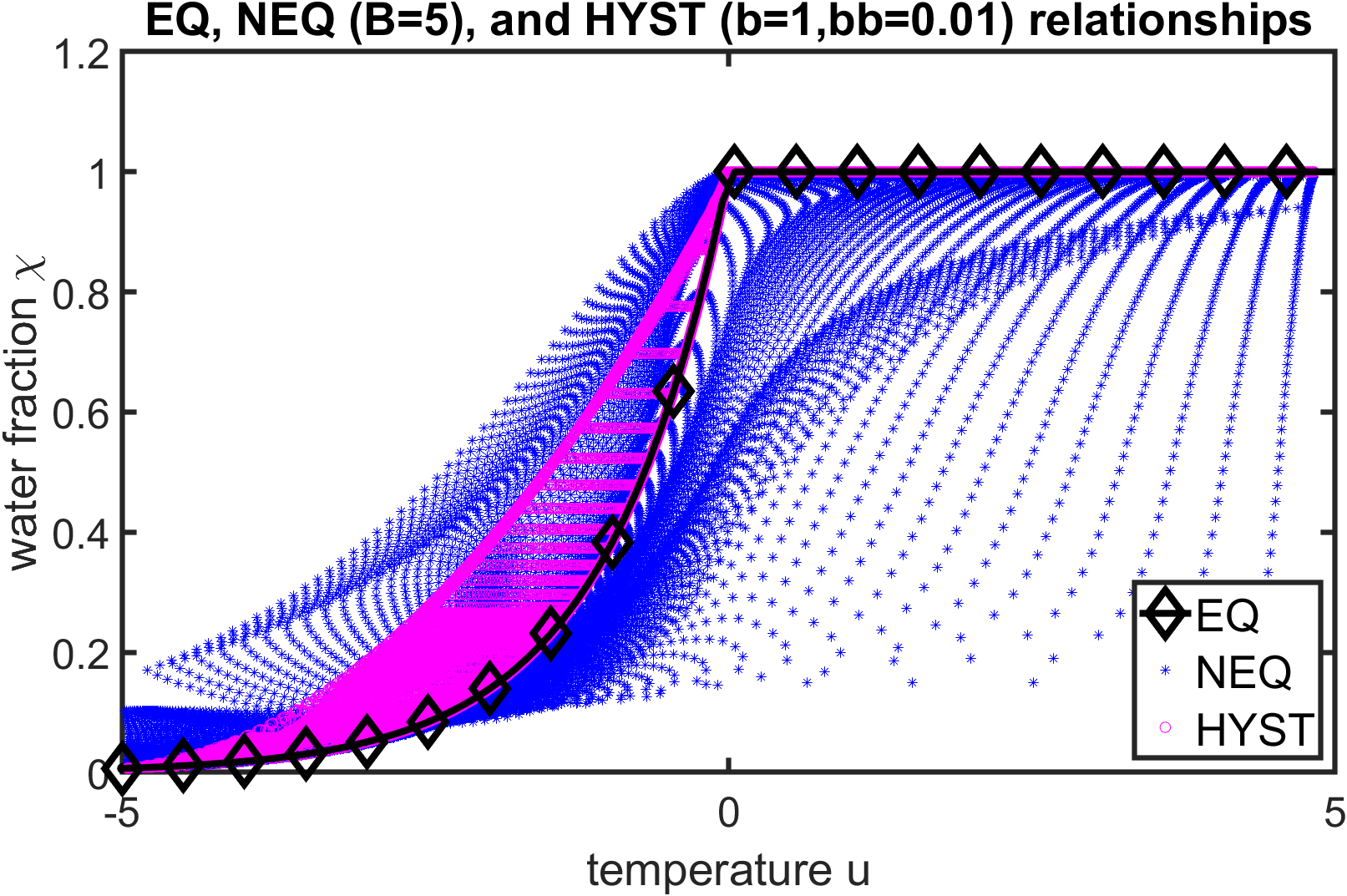}
\includegraphics[width=.45\linewidth]{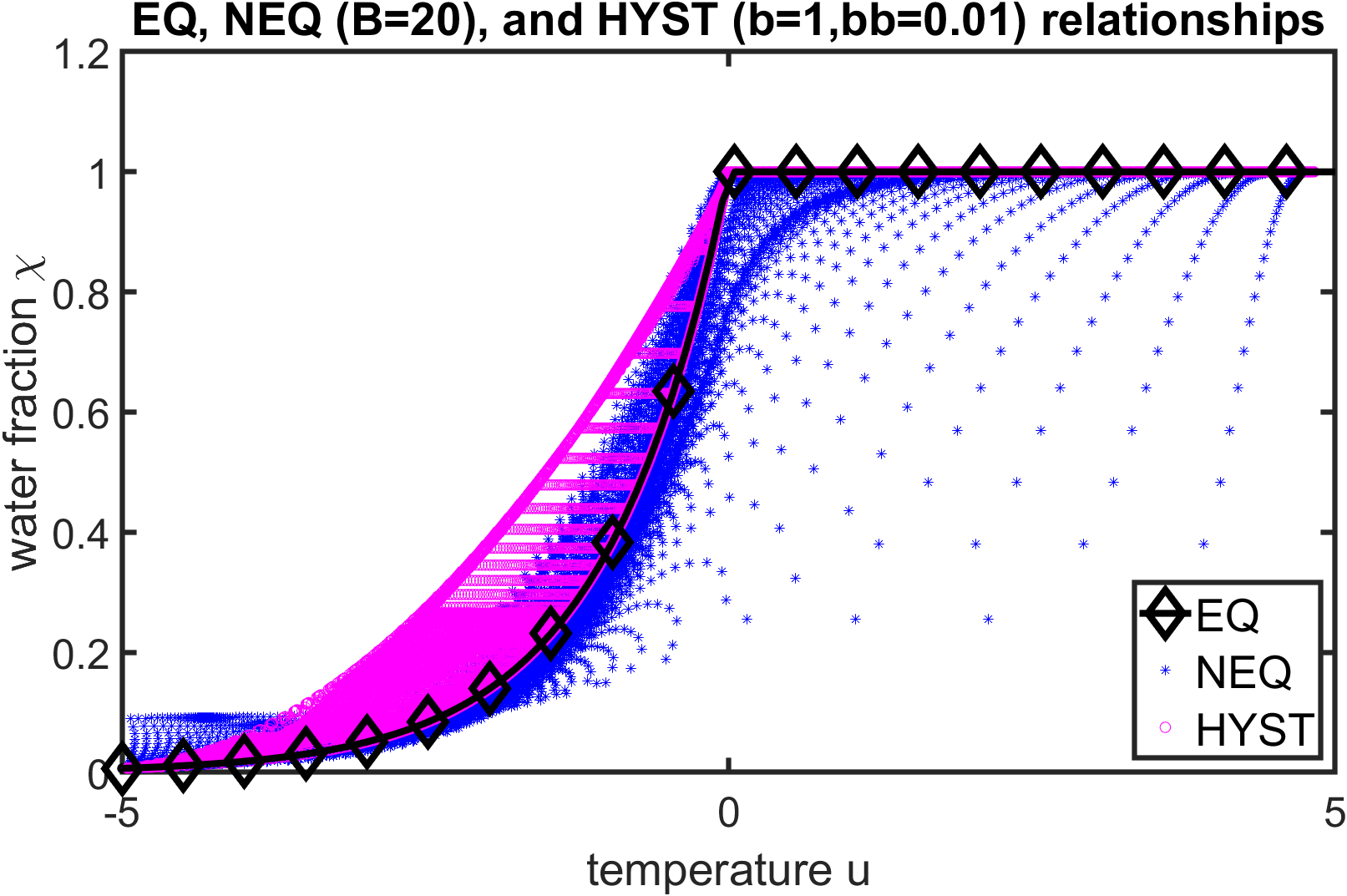}
\caption{Phase plot of the solutions to Example~\ref{ex:pde-ALL} for $B=5$ (left) and $B=20$ (right).}
\label{fig:pde-ALL}
\end{figure}

In the end, the overall dynamics of the problem is best seen in the phase plot of 
\\
$(u(x_j,t_n), \chi(x_j,t_n))_{j,n}$ given in Figure~\ref{fig:pde-ALL} and earlier in Figure~\ref{fig:uchi-ALL}. There we give the points collected from all the simulations for $(u^{EQ}(x_j,t_n), \chi^{EQ}(x_j,t_n))_{j,n}$, $(u^{NEQ}(x_j,t_n), \chi^{NEQ}(x_j,t_n))_{j,n}$, and $(u^{HYST}(x_j,t_n), \chi^{HYST}(x_j,t_n))_{j,n}$. We see the hysteresis loops for the (HYST) solutions. We also see the large difference between the (EQ) and the (NEQ) case, and a significant more spread for the (NEQ) case when $B=1$ and the small spread when $B=5$. 

%%%
\subsubsection{Solver performance}
The solver is quite robust even though the initiation of the free boundary around $t=0$ presents a challenge. Also, more iterations are needed around the time of dramatic change in the boundary conditions such as around $t=2$.  We record the number of required iterations needed to get the residual to satisfy $\norm{T^{(r)}}{\infty}\leq 10^{-8}$, with maximum $N_{max}$ of iterations in the more difficult time steps, down to minimum $N_{min}$ when conditions stabilize, with average number $N_{ave}$. We report these in Table~\ref{tab:pde}. We do not iterate with more than $20$ iterations. 

As expected, thanks to the analysis in Sections~\ref{sec:analysis}, there is no significant difference in performance between EQ, and NEQ, with the NEQ case appearing to be slightly ``easier'', due to the delay effects of relaxation. However, HYST cases seem to require more effort; we believe this is because  
they feature lower degree of semi-smoothness. 

There is some but small difference in the dependence on the dimension $M$ (primarily in $N_{init}$ for a fixed $\tau$. There is a larger difference for fixed $M$ on the time step $\tau$. Unlike in other nonlinear PDEs where smaller $\tau$ improves Newton convergence and that coarse grids requires fewer iterations, we see here that  the cases with coarse discretization but small time step seem to struggle more, especially for HYST case. We believe this is related to the fact that the solver has difficulty getting the free boundary right, and thus requires a larger number of iterations. Clearly more work is needed on the solvers.

%%%
\begin{table}
\centering
\caption{Number of iterations of nonlinear solver in Examples~\ref{ex:pde-EQ} and \ref{ex:pde-ALL}.  }
\label{tab:pde}
\begin{tabular}{|l|l|lll|lll|lll|}
\hline
&&
\multicolumn{3}{c|}{EQ}&
\multicolumn{3}{c|}{NEQ}&
\multicolumn{3}{c|}{HYST}
\\
$M$&$\tau$ 
& $N_{min}$& $N_{max}$& $N_{ave}$
& $N_{min}$& $N_{max}$& $N_{ave}$
& $N_{min}$& $N_{max}$& $N_{ave}$
\\ 
\hline
100 & 0.01 & 3 & 8 & 4.34333 & 3 & 8 & 4.28333 & 3 & 12 & 6.88333
\\
100 & 0.1 & 4 & 10 & 6.26667 & 4 & 9 & 6 & 4 & 12 & 7.9
\\
50 & 0.01 & 3 & 8 & 4.28667 & 3 & 7 & 4.22667 & 3 & 20 & 7.63667
\\
50 & 0.1 & 4 & 10 & 6.1 & 4 & 9 & 5.96667 & 4 & 12 & 7.9
\\
20 & 0.01 & 3 & 20 & 9.87333 & 3 & 16 & 4.10333 & 3 & 20 & 12.9567
\\
20 & 0.1 & 5 & 9 & 6.03333 & 4 & 9 & 5.96667 & 4 & 20 & 8.43333
\\
\hline
\end{tabular}
\end{table}

%%%%%%%%%%%%

%%%%%%%%%%%%%
\section{Summary and acknowledgements}
\label{sec:summary}
In this paper we outlined the basic ingredients of a model for phase transitions out of equilibrium for the applications to thawing and freezing in permafrost soils. Our focus was on a unified presentation and analysis of fully implicit schemes for these phenomena. 

The main challenges are the presence of two nonlinearities $F(U)$ describing the water fraction due to the temperature, as well as $A(U)U$ associated with the nonlinear dependence of the heat conductivity on the temperature $u$. The analysis of the schemes we conducted reveals that handling the non-equilibrium and hysteretic relationships can be done in a similar fashion as of those for the equilibrium model. Computational examples show all models are robust, but the HYST case requires slightly more work and is less smooth.

In this paper we presented only a few examples.  Work is underway on identifying the appropriate data for the non-equilibrium and hysteretic models depending on the spatial scale of interest. In addition, in this paper we considered fixed formulas for $c(u),k(u)$ not directly involving the variable liquid fraction subject to the lack of equilibrium. We plan to extend our work in this direction in the future.

{\bf Acknowledgments.} 
This research was partially supported by the grants NSF DMS-1912938 ``Modeling with Constraints and Phase Transitions in Porous Media'', and NSF DMS-2309682 ``Computational mathematics of Arctic processes''. \mpc{We also would like to thank the reviewers and editors whose remarks helped to improve this manuscript.}

%%%%%%%%%%%%%
%\section{Appendix}
%%%%%%%%%%%%%%%%
%\section*{References}
\bibliographystyle{plain}
\bibliography{NONEQ}

\begin{thebibliography}{10}

\bibitem{APS}
A.~Alhammali, M.~Peszynska, and C.~Shin.
\newblock {Numerical analysis of a mixed finite element approximation of a
  coupled system modeling biofilm growth in porous media with simulations}.
\newblock {\em {IJNAM}}, 21:20--64, 2024.

\bibitem{AHan2}
Kendall Atkinson and Weimin Han.
\newblock {\em Theoretical numerical analysis}, volume~39 of {\em Texts in
  Applied Mathematics}.
\newblock Springer, New York, second edition, 2005.
\newblock A functional analysis framework.

\bibitem{BPV}
Lisa Bigler, Malgorzata Peszynska, and Naren Vohra.
\newblock {Heterogeneous {Stefan} problem and permafrost models with P0-P0
  finite elements and fully implicit monolithic solver}.
\newblock {\em Electronic Research Archive}, 30(4):1477--1531, 2022.

\bibitem{Brezis73}
H.~Br{\'e}zis.
\newblock {\em Op\'erateurs maximaux monotones et semi-groupes de contractions
  dans les espaces de {H}ilbert}.
\newblock North-Holland Publishing Co., Amsterdam, 1973.
\newblock North-Holland Mathematics Studies, No. 5. Notas de Matem{\'a}tica
  (50).

\bibitem{JiangNochetto98}
Xun Jiang and Ricardo Nochetto.
\newblock A {P}1--{P}1 finite element method for a phase relaxation model {I}:
  Quasiuniform mesh.
\newblock {\em Siam Journal on Numerical Analysis - SIAM J NUMER ANAL}, 35, 06
  1998.

\bibitem{KrasPokr89}
M.~A. Krasnoselskii and A.~V. Pokrovski\u\i.
\newblock {\em Systems with hysteresis}.
\newblock Springer-Verlag, Berlin, 1989.
\newblock Translated from the Russian by Marek Niezg\'odka.

\bibitem{LittShow95}
T.~D. Little and R.~E. Showalter.
\newblock The super-{S}tefan problem.
\newblock {\em Internat. J. Engrg. Sci.}, 33(1):67--75, 1995.

\bibitem{Mayergoyz}
I.~D. Mayergoyz.
\newblock {\em Mathematical models of hysteresis}.
\newblock Springer-Verlag, New York, 1991.

\bibitem{MedinaP18}
F.~Patricia Medina and Malgorzata Peszynska.
\newblock Stability for implicit–explicit schemes for non-equilibrium kinetic
  systems in weighted spaces with symmetrization.
\newblock {\em Journal of Computational and Applied Mathematics}, 328:216--231,
  2018.

\bibitem{Michalowski93}
Radoslaw~L. Michalowski.
\newblock A constitutive model of saturated soils for frost heave simulations.
\newblock {\em Cold Regions Science and Technology}, 22(1):47--63, 1993.

\bibitem{NicolskyRomanovskyTipenko07}
Dmitry Nicolsky, Vladimir Romanovsky, and Gennadiy Tipenko.
\newblock Using in-situ temperature measurements to estimate saturated soil
  thermal properties by solving a sequence of optimization problems.
\newblock {\em The Cryosphere}, 1, 11 2007.

\bibitem{P97}
M.~Peszynska.
\newblock A differential model of adsorption hysteresis with applications to
  chromatography.
\newblock In Jorge~Guinez Angel Domingo~Rueda, editor, {\em III Coloquio sobre
  Ecuaciones Diferenciales Y Aplicaciones, May 1997}, volume~II. Universidad
  del Zulia, 1998.

\bibitem{PHV}
M.~Peszynska, Z.~Hilliard, and N~Vohra.
\newblock Coupled flow and energy models with phase change in permafrost from
  pore- to {Darcy} scale: modeling and approximation.
\newblock {\em Journal of Computational and Applied Mathematics}, 450:115964,
  November 2024.

\bibitem{PJW02}
M.~Peszynska, E.~Jenkins, and M.~F. Wheeler.
\newblock Boundary conditions for fully implicit two-phase flow model.
\newblock In Xiaobing Feng and Tim~P. Schulze, editors, {\em Recent Advances in
  Numerical Methods for Partial Differential Equations and Applications},
  volume 306 of {\em Contemporary Mathematics Series}, pages 85--106. American
  Mathematical Society, 2002.

\bibitem{PShin21}
M.~Peszynska and C.~Shin.
\newblock Stability of a numerical scheme for methane transport in hydrate zone
  under equilibrium and non-equilibrium conditions.
\newblock {\em Computational Geosciences}, 5:1855--1886, 2021.

\bibitem{PS20}
M.~Peszynska and R.~Showalter.
\newblock Approximation of scalar conservation law with hysteresis.
\newblock {\em {SIAM Journal Numerical Analysis}}, 58(2):962--987, 2020.

\bibitem{PS98}
M.~Peszynska and R.~E. Showalter.
\newblock A transport model with adsorption hysteresis.
\newblock {\em Differential Integral Equations}, 11(2):327--340, 1998.

\bibitem{PS21}
Malgorzata Peszynska and Ralph~E Showalter.
\newblock Approximation of hysteresis functional.
\newblock {\em Journal of Computational and Applied Mathematics}, 389:113356,
  2021.

\bibitem{PVB}
Malgorzata Peszynska, Naren Vohra, and Lisa Bigler.
\newblock {Upscaling an Extended Heterogeneous Stefan Problem from the
  Pore-Scale to the Darcy Scale in Permafrost}.
\newblock {\em Multiscale Modeling \& Simulation}, 22(1):436--475, 2024.

\bibitem{RomanovskyOsterkamp2000}
Vladimir Romanovsky and Tom Osterkamp.
\newblock Effects of unfrozen water on heat and mass transport in the active
  layer and permafrost.
\newblock {\em Permafrost and Periglacial Processes}, 11:219--239, 07 2000.

\bibitem{Show-monotone}
R.~E. Showalter.
\newblock {\em Monotone operators in {B}anach space and nonlinear partial
  differential equations}, volume~49 of {\em Mathematical Surveys and
  Monographs}.
\newblock American Mathematical Society, Providence, RI, 1997.

\bibitem{Ulbrich}
Michael Ulbrich.
\newblock {\em Semismooth {N}ewton methods for variational inequalities and
  constrained optimization problems in function spaces}, volume~11 of {\em
  MOS-SIAM Series on Optimization}.
\newblock Society for Industrial and Applied Mathematics (SIAM), Philadelphia,
  PA, 2011.

\bibitem{Visintin94}
Augusto Visintin.
\newblock {\em Differential models of hysteresis}, volume 111 of {\em Applied
  Mathematical Sciences}.
\newblock Springer-Verlag, Berlin, 1994.

\bibitem{Visbook}
Augusto Visintin.
\newblock {\em Models of phase transitions}, volume~28 of {\em Progress in
  Nonlinear Differential Equations and their Applications}.
\newblock Birkh\"{a}user Boston, Inc., Boston, MA, 1996.

\bibitem{VP-Tp}
Naren Vohra and Malgorzata Peszynska.
\newblock Robust conservative scheme and nonlinear solver for phase transitions
  in heterogeneous permafrost.
\newblock {\em Journal of Computational and Applied Mathematics}, 442:115719,
  2024.

\end{thebibliography}

\end{document}